\documentclass[11pt,reqno,tbtags]{amsart}

\usepackage{main_formatting}
\begin{document}

\title[Scaling limits for critical inhomogeneous random graphs]{Revisiting scaling limits for critical inhomogeneous random graphs with finite third moments} 

\author[Addario-Berry]{Louigi Addario-Berry}
\email{louigi.addario@mcgill.ca}
\address{McGill University, Montr\'eal, Canada}
\author[Bell]{Sasha Bell}\email{sasha.bell@mail.mcgill.ca}
\address{McGill University, Montr\'eal, Canada}
\author[Deka]{Prabhanka Deka} \email{prabhanka@gmail.com} \address{Peking University, Beijing, China}
\author[Donderwinkel]{Serte Donderwinkel}
\email{s.a.donderwinkel@rug.nl}
\address{University of Groningen, Groningen, The Netherlands}
\author[Maniyar]{Sourish Maniyar}
\email{sourishparag@iisc.ac.in}
\address{Indian Institute of Science, Bengaluru, India}
\author[Wang]{Minmin Wang}
\email{minmin.wang@sussex.ac.uk}
\address{University of Sussex, Brighton, UK}
\author[Winter]{Anita Winter}
\email{anita.winter@uni-due.de}
\address{Universit\"at Duisburg-Essen, Germany}


\begin{abstract} 
We consider the rank-1\ inhomogeneous random graph in the Brownian regime in the critical window. Aldous studied the weights of the components, and showed that this ordered sequence converges in the $\ell^2$-topology to the ordered excursions of a Brownian motion with parabolic drift when appropriately rescaled \cite{Aldous97}, as the number of vertices $n$ tends to infinity. We show that, under the finite third moment condition, the same conclusion holds for the ordered component sizes. This in particular proves a result claimed by Bhamidi, Van der Hofstad and Van Leeuwaarden \cite{BHL2010}. 
We also show that, for the large components, the ranking by component weights coincides with the ranking by component sizes with high probability as $n \to \infty$.  
\end{abstract}

\maketitle


\section{Introduction and main result}\label{intro}

This work concerns the {\em rank-1 inhomogeneous random graph}. This random graph model was introduced by Aldous in \cite{Aldous97} and has since been heavily studied   \cite{MR4911757, BHL2010,MR4586227,MR4683376,MR4474511,broutin2021limits,BSW17,MR3622886,MR3068034,BHL12,MR2337396,MR2213964} under various names, including the {\em Poissonian graph}, the {\em Norros--Reittu model}, and {\em the multiplicative random graph}. 
It may be defined as follows. A {\em weight sequence}\label{def:wt seq} is a vector $\mathrm{w}\in (0,\infty)^n$ for some $n\ge 1$. For each $n\geq 1$, fix a weight sequence $\wtv \coloneqq(\wt{1},\ldots,\wt{n})$, and for $t_n>0$, construct a random graph with vertex set $[n]\coloneqq\{1,\ldots,n\}$ by including an edge between vertices $u$ and $v$ with probability
\begin{equation}\label{eq:puvdef}
p_{uv}= p_{u,v}(\wtv; t_n)\coloneqq 1-\exp\left\{-\frac{t_n\wt{u}\wt{v}}{\sum_{i\in [n]} \wt{i}}\right\} 
\end{equation} 
 independently for each distinct pair $u,v$ of elements of $[n]$.

We fix $\cw\in \R$, and set $t_n=1+\cw n^{-1/3}$ for $n$ large enough that $t_n$ is positive. We call the random graph that we thus obtain $\Wn$.\label{def:rank one graph} Under suitable assumptions on the weight sequences $\wtv$, the graph $\Wn$ then falls within the so-called {\em critical regime}, characterised by first-order randomness in the sizes of the largest connected components.  
The most well-known special case is that of the {\em Erd\H{o}s-R\'enyi random graph}, obtained by taking $\wt{i}=1$ for all $i \in [n]$. In this case, it has been known since the 1990s that the largest connected components of $\Wn$ have size of order $n^{2/3}$ \cite{MR1220220,Aldous97}.

\subsection{Our contributions} 
 In \cite{Aldous97}, Aldous considered the rank-1 model at criticality in the Brownian regime, but instead of counting vertices in a connected component, he studied the sum of the weights of the vertices in the component, and proved a Brownian-type limit theorem for the ranked sequence of the component weights.  
This result raises the following questions, both of which we answer affirmatively in our work, in the finite third moment setting:
\begin{enumerate}
\item
(Distributional ranking consistency) If we substitute component weights by component sizes, does the same limit theorem hold?
\item
(Almost sure ranking consistency) Is it the case that, for the large components, the ranking by component weight coincides with the ranking by size, with high probability as the graph size tends to infinity?
\end{enumerate}
We shall assume the following conditions on the progression of weight sequences $(\wtv,n\ge 1)$. Here and below, for $\mathrm{w}=(w_1,\ldots,w_n) \in (0,\infty)^n$ and $c>0$, we write $|\mathrm{w}|_c=(\sum_{i \in [n]}|w_i|^c)^{1/c}$.
\begin{enumerate}
\item\textbf{Maximal weight: } As $n\to \infty$,
\begin{equation}\label{maximal weight} 
    \max_{i\in[n]} \wt{i} =  o(n^{1/3}).
    \end{equation}
\item\label{moment convergence} \textbf{Convergence of moments:}  There exist $\mi,\miii\in(0,\infty)$ such that as $n\to \infty$,
\begin{align}
&\frac{1}{n}|\wtv|_1 = \mi+o(n^{-1/3}), \label{first moment} \\
&\frac{1}{n}|\wtv|_2^2 = \mi+o(n^{-1/3}), \label{second moment}\\ 
&\frac{1}{n}|\wtv|_3^3 = \miii+o(1). \label{third moment}
\end{align}
\end{enumerate}
The choice of $t_n=1+\lambda n^{-1/3}$ corresponds to the critical window around the point $t_n=1$ where the size of the largest component undergoes a phase transition. It was shown in \cite{MR2337396} that,  if $t_n=t$, not depending on $n$, the size of the largest component has a deterministic limit in probability when rescaled by $n$, which is positive if and only if $t>1$, and that the size of the remaining components vanishes when rescaled by $n$. We show that at $t=1$, for any $k\ge 1$, the $k$th largest component has size and weight of order $n^{2/3}$ and has a non-degenerate limit on this scale, and that this is also true for $t_n=1+\lambda n^{-1/3}$ with $\lambda\in \R$, albeit with a different random limit. Moreover, our result implies that, on the scale $n^{2/3}$, the size of the largest component diverges if $\lambda=\lambda_n\to\infty$ and vanishes if $\lambda=\lambda_n \to - \infty$, so that this range of values of $t_n$ indeed covers the full window where the behaviour at the critical point $t_n=1$ persists.

More generally, we can consider the case where, instead of~\eqref{second moment}, we have
$\frac{1}{n}|\wtv|_2^2 = \mi''+o(n^{-1/3})$ for some $\mi''\in (0, \infty)$. Replacing 
$t_n$ by $t_n^{''}=\tfrac{\mu}{\mu''}(1+\lambda n^{-1/3})$, we obtain a graph with the same law as $\mathcal{G}^\cw(\tfrac{\mu}{\mu''}\wtv)$. Thus, up to a renormalisation of the weights, we can always assume that the critical point is at $t_n=1$.

We will now formally state our results. For a connected component $C$ of $\Wn$, the {\em size} of $C$, denoted $|C|$, is the number of vertices of $C$.\label{def:size of C} The {\em weight} of $C$, denoted $w(C)$, is $\sum_{v \in V(C)} \wt{v}$.\label{def:wt of C} Our main result is the following distributional limit theorem. 
Denote by $\ell^2=\{\mathbf x=(x_i)_{i\in \N}\in \R^\N: \sum_ix_i^2<\infty\}$ the space of square-summable sequences equipped with the $\ell^2$-norm: $|x|_2=(\sum_{i\in \N}x_i^2)^{1/2}$. For $d\in \N$, we then define $\ell^2(\R^d)=(\ell^2)^d$. 

\begin{thm}\label{samelaw} 
Suppose $(\wtv,n \ge 1)$ satisfies \labelcref{maximal weight}--\labelcref{third moment}. For $n\ge 1$, list the components of $\Wn$ in decreasing order of size as $(\lgcl{i},i \ge 1)$, and in decreasing order of  weight  as $(\hvcl{i},i \ge 1)$, in both cases breaking ties by smallest vertex label in the component.
Then there exists a random $\ell^2(\R)$-valued sequence $(\excsize{i},i \ge 1)$ such that 
\begin{equation}
\label{eqcv: main}
\big(n^{-2/3}\big(|\lgcl{i}|,|\hvcl{i}|,w(\lgcl{i}),w(\hvcl{i})\big),i\geq 1\big)  \to ((\excsize{i},\excsize{i},\excsize{i},\excsize{i}),i\geq 1)
\end{equation}
in distribution in $\ell^2(\R^4)$, as $n \to \infty$.
\end{thm}
Note that \eqref{first moment}, \eqref{second moment} and the Cauchy--Schwarz inequality imply that $\mu\le 1$. We observe that if $\mu<1$, then the sum of the component weights is asymptotically smaller than the sum of the component sizes. Thus, for $\mu < 1$, the asymptotic equivalence of the component weights and sizes clearly cannot hold for \emph{all} the components of $\Wn$. That this equivalence holds for the largest and heaviest components is a consequence of the criticality condition. 

The work \cite{BHL2010} claims to show the convergence of $(n^{-2/3}|\lgcl{i}|,i\ge 1)$, under an additional regularity condition on $\wtv$, but that work contains an error in applying Aldous' limit theory of size-biased point processes; we elaborate on this in Section~\ref{sec: l2}, just before the proof of Theorem~\ref{samelaw}. 
The work~\cite{Turova}, which considers a setting similar to~\cite{BHL2010}, contains the same type of mistake.
Furthermore, the work \cite{BSW17} relies on the results in \cite{BHL2010}, so with our new proof we also show that these results indeed hold.

It will follow from an explicit description of the limiting sequence $(\excsize{i},i \ge 1)$ that almost surely $\excsize{i+1}<\excsize{i}$ for all $i$. Using this fact, the following corollary of Theorem~\ref{samelaw} is immediate.
\begin{cor}\label{cor_consist}
For all $i\ge 1$, with high probability, $\hvcl{i} = \lgcl{i}$.
\end{cor}
In the corollary, and hereafter, the phrase ``with high probability" means ``with probability tending to $1$ as $n$ tends to infinity".  Observe that to derive the corollary, it is sufficient that \eqref{eqcv: main} holds in the product topology. This is important because we will first show this weaker result and then use \cref{cor_consist} for the proof of \cref{samelaw}.

We now give a slightly informal description of the limit object $(\excsize{i},i \ge 1)$; a more precise definition appears at the start of \cref{S:proof}. The law of $(\excsize{i},i \ge 1)$ depends on $\lambda,\mu$ and $\mu'$ as in \labelcref{eq:puvdef}-\labelcref{third moment}, but we suppress this dependence in the notation. 

Let $(\sbm(t),t\geq 0)$\label{def:sbm} be a standard Brownian motion, and for $t\ge 0$ let
\begin{equation}
\label{defW}
\bm(t) = \bm(t; \cw, \mu,\mu')\coloneqq \sqrt{\frac{\mu'}{\mu}}\sbm(t) + \cw t - \frac12 \frac{\mu'}{\mu^2}t^2,
\end{equation}
so that $(\bm(t),t\ge 0)$ is an inhomogeneous Brownian motion with drift $\lambda-t\mu'/\mu^2$ at time $t$. Then, for $i\geq 1$, let $\excsize{i}$\label{def:excsize} be the length of the $i$th largest excursion of $\bm(t)$ above its running infimum.

\subsection{Relation to prior work}\label{sec:relation}

In \cite{Aldous97}, Aldous shows that, under the same conditions as used in this work, the component weights converge when appropriately rescaled. Aldous' work was later generalised to the L\'evy setting by Aldous and Limic~\cite{AlLi98}.
In~\cite{BHL2010}, the authors investigated distributional ranking consistency under the assumptions~\eqref{maximal weight}-\eqref{third moment} on the vertex weights, $t_n=1+\lambda n^{-1/3}$ and an additional regularity condition on the weight sequence; see also~\cite{Turova}. (As mentioned, the proof in \cite{BHL2010} is incorrect; see Remark~\ref{rem: SBPP} in Section~\ref{sec: l2}.) Bhamidi, van der Hofstad and van Leeuwaarden~\cite{BHL12} considered the same question under the assumption that the vertex weights follow a specific power-law tail, and gave an affirmative answer in that case. Broutin, Duquesne and the sixth author~\cite{broutin2021limits} studied the model under the general asymptotic regime considered in Aldous and Limic~\cite{AlLi98} and proposed a sufficient condition for a positive answer to the first question with respect to the product topology (see Theorem 2.8 and Remark 2.2 there). We also point out that the arguments in~\cite{broutin2021limits} actually prove almost sure consistency under this condition.

We next give a brief account of our approach. 
Our analysis relies on two key ingredients. First, we introduce a coupling of two breadth-first exploration processes, one which visits components in a size-biased order and the other which visits components in a weight-biased order. In~\cref{sizes conv T} we show that, under the coupling, the sizes and weights of the heaviest components discovered amongst the first $\Theta(n^{2/3})$ vertices in the two explorations have an identical scaling limit. Second, \cref{rhino} shows that in the size-biased exploration we are unlikely to have missed a large component by time $\Theta(n^{2/3})$ and in the weight-biased exploration we are unlikely to have missed a heavy component by time $\Theta(n^{2/3})$. Combining these two ingredients with some careful analysis then allows us to prove distributional ranking consistency in the product topology and almost sure ranking consistency.  Let us mention that similar results to \cref{rhino} have been obtained in the literature using different, and sometimes model-dependent methods, including Aldous’ size-biased point process technique~\cite{Aldous97}, and a large deviation approach in~\cite{BHL12}; see also the susceptibility argument for the configuration model by Janson~\cite{Jan10}. Our proof of \cref{rhino} is elementary, and we believe it can be adapted to more general situations. 

To improve the convergence to the $\ell^2$-topology, we propose a generalisation of Aldous’ size-biased point process technique~\cite{Aldous97}. In its original form, the technique, which Aldous used to extract $\ell^2$-convergence of the component weights from the convergence of exploration processes, requires that the excursion lengths correspond to the component weights (resp.~sizes) and the components are explored in a weight-biased (resp.~size-biased) ordering. In many examples, however, both requirements are not fulfilled. For instance, in~\cite{BHL2010}, although the excursion lengths encoded the number of vertices, the exploration was conducted in a weight-biased order. In such a setting, the application of the size-biased point process technique in~\cite{BHL2010} (used in their proof for product topology as well as $\ell^2$-topology convergence) seems complicated. 
In our size-biased exploration introduced in Section~\ref{defining St}, we use the excursion lengths to encode the weights of the components that appear in a size-biased fashion. The main difference in our approach comes from the fact that with the help of an auxiliary counting process, we then can construct (in Section~\ref{sec: l2}) the ``correct'' size-biased point process. Once again, we believe that combined with the coupling of two exploration processes, this version of the size-biased point process technique can be applied to other regimes for the vertex weights as well. In fact, we conjecture the following.

{\bf Conjecture.} Distributional ranking consistency and almost sure ranking consistency also hold under the assumptions of~\cite[Proposition 7]{AlLi98}. To be precise, we conjecture that under those assumptions, the convergence of the component weights in $\ell^2$ topology stated in the proposition also holds for the sequence of component sizes and, moreover, that with high probability the $k$th largest connected component is also the $k$th heaviest as $n\to\infty$. 

\subsection{Weight- and size-biased exploration processes of $\Wn$}
\label{defining St}
In this section, we describe how to construct a graph $G_n$\label{def:Gn wt} distributed as $\Wn$ based on
the following ingredients:
\begin{itemize}
\item a root-selector map $r:\,2^{[n]}\setminus\{\varnothing\}\to[n]$ such that $r(V)\in V$ for all $V\subseteq[n]$,
\item an independent family $\{\jump_v:\,v\in[n]\}$ with $\jump_v\disteq\mathrm{Exp}\Big(\frac{\wt{v}(1+\cw n^{-1/3})}{|\wtv|_1}\Big)$, and
\item an independent family $\{\jump_{v,v'}:\,v\not = v'\in[n]\}$ with $\jump_{v,v'}\disteq\mathrm{Exp}\Big(\frac{\wt{v}\wt{v'}(1+\cw n^{-1/3})}{|\wtv|_1}\Big)$. 
\end{itemize}
The construction is essentially a breadth-first search exploration of the connected components of the graph, which uses the root-selector map for deciding which vertex to explore from when starting the exploration of a new component; the random variables $\{J_v,v \in [n]\}$ for constructing the breadth-first search trees of the components; and the random variables $\{J_{v,v'}:v \ne v' \in [n]\}$ for determining the remaining edges of the graph. Such explorations have seen substantial use in the study of random graphs and of related coalescent processes \cite{Aldous97,AlLi98,BHL2010,BHL12,BSW17,broutin2021limits,MR4474511,MR3622886,limiceternal,Turova,MR3068034}.

We begin by declaring $r_1:=r([n])$ as the root of the first component. 
Let $\mathrm{ch}(r_1):=\{v\in [n]\setminus \{r_1\},J_v<w^{(n)}_{r_1}\}$ be the {\em set of children} of $r_1$. Put $b(1):=r_1$, and write $(b(2),\ldots,b(1+|\mathrm{ch}(r_1)|))$ for the list of elements of $\mathrm{ch}(r_1)$, ordered so that $J_{b(2)}<\cdots<J_{b(|\mathrm{ch}(r_1)|+1)}$. We 
next let $\mathrm{ch}(b(2)):=\{v\in [n]\setminus\{r_1\}:\,w^{(n)}_{r_1}\le J_v<w^{(n)}_{r_1}+w^{(n)}_{b(2)}\}$ be the set of the children of $b(2)$, and write $(b(2+|\mathrm{ch}(r_1)|),...,b(1+|\mathrm{ch}(r_1)| + |\mathrm{ch}(b_2)|))$ for the list of elements of $\mathrm{ch}(b(2))$ ordered so that $J_{2+b(|\mathrm{ch}(r_1)|)}<...<J_{b(1+|\mathrm{ch}(r_1)| + |\mathrm{ch}(b_2)|)}$. Until we can not find more descendants of $r_1$, we continue to explore the sets of children  
\begin{equation}
\label{e:children}
    \mathrm{ch}(b(k)):=\Big\{v\in [n]\setminus\{r_1\}:\,\sum_{l\in[k-1]}w^{(n)}_{b(l)}\le J_{v}<\sum_{l\in[k]}w^{(n)}_{b(l)}\Big\}
\end{equation}
of each of the vertices $b(k)$ we have discovered and in the order we have discovered them, always ordering the children of a given vertex in increasing order of the $J_v$ values.

We write $V^r_1$ for the set consisting of $r_1$ and all its descendants, and refer to it as the first component of $G^r_n$. 
We assign each vertex of $V^r_1$ a {\em discovery time}, as follows. We say that vertex $r_1=b(1)$ was discovered at time $0$, and for $1 < k \le |V_1^r|$, that vertex $b(k)$ was discovered at time $J_{b(k)}$. Also, for $1 \le k \le |V_1^r|$, we say that vertex $b(k)$ was {\em explored} during the time interval 
\begin{equation}\label{eq:explore_interval}
\Big[\sum_{l \in [k-1]} \wt{b(l)},\sum_{l \in [k]} \wt{b(l)}\Big)\, ,
\end{equation}
which has length $w_{b(k)}$.
These intervals are pairwise disjoint for distinct elements of~$V_1^r$.

If $V^r_1 \ne [n]$ then we declare $r_2=r([n]\setminus V^r_1)$ to be the root of the second component and set $b(|V_1^r|+1):=r_2$. Once more we use breadth-first search on the remaining vertices $[n]\setminus V^r_1$ 
to construct a tree rooted at $r_2$ together with a breadth-first ordering of its vertices. Refer to the set $V^r_2$ consisting of
$r_2$ and all its descendants as the second component of $G^r_n$; write $(b(|V_1^r| + 1),\ldots,b(|V_1^r|+|V_2^r|))$ for the list of its elements in breadth-first order. We say that vertex $r_2=b(|V_1^r|+1)$ was discovered at time $\sum_{l \in V_1^r} \wt{b(l)}$, and for $|V_1^r|+1<k\le |V_1^r \cup V_2^r|$ say that vertex $b(k)$ was discovered at time $J_{b(k)}$. Also, for $|V_1^r|<k\le |V_1^r \cup V_2^r|$ again say that $b(k)$ was explored during the time interval \eqref{eq:explore_interval}. Proceed in this manner until there are no more unexplored vertices, and let $b:[n]\to[n]$ be the resulting order of the vertices of $G_n^r$. Denote by 
$K^r_n$ the total number of components of the resulting graph.  

We use the family $\{\jump_{v,v'}:\,v\not= v'\in[n]\}$ to extend the obtained forest to a graph by adding an edge between any two vertices $v,u\in V$ with $J_{v,u}<1$ if $\jump_u<\jump_v$ and $\jump_v\le\sum_{l=1}^{b^{-1}(u)-1}w^{(n)}_{b(l)}$, when $u$ is due to be  explored. 
This yields a random graph, $G^r_n$, which is distributed as $\Wn$; see \cite[Section 3.1]{Aldous97}, \cite[Section 2]{limiceternal}. 

Importantly, 
the specific choice of the root selector map does not affect the 
distribution of the resulting random graph. 
This implies that we can even choose the root selector map in a random fashion, 
provided that it is independent of the families 
$\{\jump_v:\,v\in[n]\}$ and $\{\jump_{v,v'}:\,v\not = v'\in[n]\}$, 
and the resulting graph is still $\Wn$-distributed. 
Moreover, by varying the choice of root selector map, while using the same random variables $\{J_v:\,v \in [n]\}$ and $\{J_{v,v'}:\,v \ne v' \in [n]\}$, we can construct different exploration processes which construct $\Wn$-distributed graphs that are coupled to one another. This fact undergirds much of our proof; as alluded to in Section~\ref{sec:relation}, we will use two root selector maps, which yield versions of the exploration that explore the components in size-biased order and in weight-biased order.

Since our focus is on  the weights and the sizes of the components of $\Wn$, rather than on the graph-theoretic structure of $\Wn$, it is sufficient for our purposes to give a representation of the vertex sets of the connected components of the breadth-first spanning forest of $\Wn$ described above. For this, consider the real-valued stochastic process $S_n=(S_n(t),t \ge 0)$ defined by
\begin{equation}\label{Sn}
\newrw(t) \coloneqq -t + \sum_{v \in [n]} \wt{v}\I{\jump_v\leq t},\quad t\geq 0,
\end{equation}
We next describe how we can use the root selector map $r$ and the stochastic process $S_n$ to recover these vertex sets. Recall that for all non-root vertices $v$, the random variables $\jump_v$ correspond to discovery times, and that each vertex $v$ is explored starting at some point after it is discovered, over a time interval of length $\wt{v}$, which determines the set of children of $v$ in the spanning forest. Therefore, 
\[
   \wst^r_1:=\inf\big\{t \ge 0:\, \wt{r_1}\I{\jump_{r_1}> t} +S_n(t)= 0\big\}
\]
is the time required to complete the exploration of 
\[
   \wvx^r_1 
   = 
   \{r_1\}\cup \big\{u \in [n]:\, J_u \in [0,\wst^r_1)\}\, , 
\]
the set of vertices of the first component. 
More generally, we can define a sequence of increasing random times 
$0=\tau^r_0<\tau^r_1<\tau^r_2<\cdots<K^r_n$ as 
\begin{equation}
\label{T1 def new}
  \wst^r_k := 
  \inf\Big\{t \ge 0:\, \sum_{l\in[k]} w^{(n)}_{r_l}\I{\jump_{r_l}> t} +S_n(t)= 0\Big\},
\end{equation}
and find that for all $k\in[K_n^r]$,
\[
\label{tauVk}
   \wvx^r_k 
   = 
   \{r_k\}\cup \big\{u \in [n]: J_u \in [\wst^r_{k-1},\wst^r_k)\}\setminus\{r_1,\ldots,r_{k-1}\big\}.
\]
Thus, we can recover the vertex sets of the components from $S$ and $r$. Also, for all $k \in [K_n^r]$, the weight $w(V^r_{k})$ can be recovered as $\wst^r_{k}-\wst^r_{k-1}$.

The proof of our main theorem relies on two  choices of random root selector map. First, 
in the {\em weight-biased} exploration process we select the root in a weight-biased fashion, i.e., for a non-empty $V\subseteq [n]$, we use the random selector $R$ with 
\[
   \Pp(R(V)=v)=\frac{w^{(n)}_v}{\sum_{v'\in V}w^{(n)}_{v'}}\I{V}(v).
\]   
In this case we abbreviate $\tau_k:=\tau^R_k$,  $V_k:=V^R_k$,  $K_n:=K_n^R$ and $G_n:=G_n^R$, write $R_1,\ldots,R_{K_n}$\label{def:wroot} for the roots of the trees in the resulting forest, listed in their order of discovery, and note for future use that in this case we may rewrite the formula for $\wst_k$ as 
\begin{equation}\label{eq:tauk}
\wst_{k}=\inf\left\{t\geq 0:\newrw(t) =-\sum_{l\in [k]} \wt{\wroot_l}\I{J_{\wroot_l} > t}\right\}
\, .
\end{equation} \label{tauV1}
For $k \ge 1$, we let $\hvclw{k}$\label{def:hvclw} and $\lgclw{k}$\label{def:lgclw} denote the vertex sets of the $k$th heaviest and $k$th largest components, respectively, of $G_n$; if $k \ge K_n$ then $\hvclw{k}=\lgclw{k}\coloneqq\varnothing$.  

It is shown in \cite[Section 3.1]{Aldous97} that 
the rooted forest we obtain is distributed as $\Wn$. (The process considered in that paper is described slightly differently, but is distributionally identical to ours. 
For more details, see \cref{conn to Aldous}.) Moreover, as is also noted in \cite{Aldous97}, since the vertices $\wroot_1,\ldots,\wroot_{\wK_n}$ are chosen in a weight-biased fashion, the resulting ordering $(\wvx_1,\ldots,\wvx_{\wK_n})$ is a weight-biased permutation of the vertex sets of the components of $\Wn$. 

Second, in the {\em size-biased} exploration process we select the root uniformly, i.e., for a non-empty set $V\subseteq [n]$, we use the random selector $\hat{R}$ with 
\[
   \Pp\big(\hat{R}(V)=v\big)=\frac{1}{|V|}\I{V}(v).
\]
This name is justified by the fact that the vertex sets $(V_i^{\hat{R}},i \in [K_n^{\hat{R}}])$ of the connected components of the resulting graph $G_n^{\hat{R}}$ appear in size-biased random order; see \cite{limiceternal}. We abbreviate $\hat{\tau}_k:=\tau^{\hat{R}}_k$,\label{def:tau k hat} $\hat{V}_k:=V^{\hat{R}}_k$,\label{def:Vk hat}   $\hat{K}_n:=K_n^{\hat{R}}$ and $\hat{G}_n:=G_n^{\hat{R}}$.\label{def:Gn hat} 
We also write $\hat{R}_1,\ldots,\hat{R}_{\hat{K}_n}$\label{def:sroot} for the roots of the trees in the resulting forest, listed in their order of discovery, and record that $\hat{\wst}_k$ is then given by 
\begin{equation}\label{eq:taukhat}
\hat{\wst}_{k}=\inf\left\{t\geq 0:\newrw(t) =-\sum_{l\in [k]} \wt{\hat{\wroot}_l}\I{J_{\hat{\wroot}_l} > t}\right\}
\, ,
\end{equation}
We let $\hvcls{k}$\label{def:hvcls} and $\lgcls{k}$\label{def:lgcls} denote the vertex sets of the $k$th heaviest and $k$th largest component, respectively, of $\hat{G}_n$, with $\hvcls{k}=\lgcls{k}\coloneqq\varnothing$ if $k \ge \hat{K}_n$.

It follows from the above discussion that both $G_n$ and $\hat{G}_n$ are $\Wn$-distributed, and thus that $\hvcls{k} \eqdist \hvclw{k}$ and $\lgcls{k} \eqdist \lgclw{k}$, but it does not need to be the case that $\hvcls{k} \eqas \hvclw{k}$ or $\lgcls{k} \eqas \lgclw{k}$.

\subsection{Connection to Aldous' process} 
\label{conn to Aldous}
In this section, we explain the connection between the breadth-first construction described in \cref{defining St} and the one given in \cite[Section 3.1]{Aldous97}. In the notation of~\cite{Aldous97}, Aldous constructs a random graph $G'_n$ on the vertex set $[n]$ using a vector $\xtv=(\xt{i})_{i\in [n]}$ and a real number $\q>0$, 
in which edges are independently present, and edge $ij$ is present with probability $1-\exp{(-\q\xt{i}\xt{j})}$. 
He then studies the asymptotic behaviour of the graphs under the hypotheses of \cite[Proposition 4]{Aldous97} on the sequence $((\xtv, \q), n\ge 1)$. 

The weight vector sequence $(\wtv, n\ge 1)$ in \cref{samelaw} can be transformed to satisfy the hypotheses of \cite[Proposition 4]{Aldous97} and, by choosing $\q$ appropriately, in such a way that the edge connection probabilities are also equal. 
Indeed, taking
\begin{equation}\label{def: an}
a_n = \frac{(\miii)^{1/3}}{\mi n^{2/3}}, \quad \xtv=a_n\wtv\text{ and }\quad \q = \frac{1}{a_n^{2}}\cdot\frac{1+\lambda n^{-1/3}}{|\wtv|_1},
\end{equation}
we find that 
\[
1-\exp(-\q \xt{i}\xt{j}) =1-\exp\left\{-\frac{(1+\cw n^{-1/3})\wt{i}\wt{j}}{|\wtv|_1}\right\},  
\]
for all $1\le i<j\le n$. 
Furthermore, if $(\wtv, n\ge 1)$ is prescribed as in~\cref{samelaw}, then the sequence $((a_n\wtv, \q), n\ge 1)$ satisfy the hypotheses of \cite[Proposition 4]{Aldous97} with the parameter $t$ in Aldous' notation given by $(\mu')^{-2/3}\mu\lambda$.

We also note that there is a minor difference between the way the breadth-first construction is described in this work and in \cite{Aldous97}. In Section~\ref{conn to Aldous}, we use exponentials $(\jump_v,v \in [n])$, and the children of a vertex $u$ 
explored starting at time $t$ are given by the set $\{v \in [n]: t \le J_v \le t+\wt{u}\}$. On the other hand, in \cite[Section 3.1]{Aldous97}, the exploration is driven by a family of exponentials $(\jump_{u,v},u,v \in [n])$, where for fixed $v$, the random variables $(\jump_{u,v},u \in [n])$ are identically distributed and have the same law as our $\jump_v$. The children of vertex $u$ are then given by the set $\{v \in [n]: \jump_{u,v} \le \wt{u}\}$. The laws of the resulting graphs and of the process $\newrw$ associated to the exploration are the same for our construction and for the construction in \cite{Aldous97}, due to the memoryless property of exponentials.

\subsection{Acknowledgements}
This project was initiated during the 2024~BIRS-CMI meeting ``Mathematical Foundations of Network Models and their Applications'' (24w4004), supported by BIRS, CMI and NHBM. LAB acknowledges the support of NSERC and of the Canada Research Chairs program. SD acknowledges the financial support of the CogniGron research center
and the Ubbo Emmius Funds (University of Groningen). Her research was
also supported by the Marie Skłodowska-Curie grant GraPhTra
(Universality in phase transitions in random graphs), grant agreement ID
101211705.
SB was supported by NSERC CGS-M and FRQNT Bourse de formation à la maîtrise.
MW's research is partially supported by the UKRI-EPSRC grant EP/W033585/1. 
AW was supported by the Deutsche Forschungsgemeinschaft (through Project-ID 444091549 within SPP-2265).
 Part of this research was performed while LAB, SB and SD were
visiting the Simons Laufer Mathematical Sciences Institute (SLMath), which is
supported by the National Science Foundation (Grant No. DMS-1928930).

\subsection{Notation}   
For sequences of random variables $(f(n),n\geq 1)$ and $(g(n),n\geq 1)$, we write $f(n) = o_\Pp(g(n))$ to indicate that $\frac{f(n)}{g(n)}\I{g(n)\ne 0} \to 0$ in probability as $n\rightarrow\infty$, and we write  $f(n) = O_\Pp(g(n))$ to indicate that the sequence $\left(\frac{f(n)}{g(n)}\I{g(n)\ne 0},n\geq 1\right)$ is tight. For a random variable $X$, we write $f(n) \convdistn X$ (resp.\ $f(n)\convprobn X$) to indicate that $f(n)$ converges to $X$ in distribution (resp.\ in probability) as $n\to \infty$. For two random variables $A$ and $B$, we write $A \eqdist B$ to indicate that $A$ and $B$ are equal in distribution.
\newcommand{\SplitTableNotations}{\end{tabular}
\end{table}

\begin{table}[htbp]
\centering
\begin{tabular}{ |c | c  | c | }
\hline
\textbf{Symbol} & \textbf{Explanation} &  \\
\hline }

\begin{center}

\begin{longtable}{@{}|p{0.17\textwidth}|>{\raggedright}p{0.62\textwidth}|p{0.05\textwidth}|@{}}

\hline \multicolumn{1}{|c|}{\textbf{Symbol}} & \multicolumn{1}{c|}{\textbf{Explanation}} & \multicolumn{1}{c|}{\textbf{Location}} \\ \hline 
\endfirsthead

\multicolumn{3}{c}
{{\bfseries  Table of notation -- continued from previous page}} \\
\hline \multicolumn{1}{|c|}{\textbf{Symbol}} & \multicolumn{1}{c|}{\textbf{Explanation}} & \multicolumn{1}{c|}{\textbf{Location}} \\ \hline 
\endhead

\hline \multicolumn{3}{|r|}{{Continued on next page}} \\ \hline
\endfoot

\hline
\endlastfoot

\hline
$\wtv$ & $\wtv=(w^{(n)}_1,\ldots,w^{(n)}_n)$, a weight sequence & Page~\pageref{def:wt seq} \\
$\Wn$ & rank-1 inhomogeneous random graph with weight sequence $\wtv$ & Page~\pageref{def:rank one graph} \\ 
$|C|$ & the size (i.e.\ the number of vertices) of a component $C$ & Page~\pageref{def:size of C}\\
$w(C)$ & the weight (i.e.\ the total weight of the vertices) of a component $C$ & Page~\pageref{def:wt of C}\\

\hline
$G_n$ & $\Wn$-distributed graph, constructed via the \em{weight-biased exploration process} & Page~\pageref{def:Gn wt} \\
$\wvx_i$ & vertex set of $i$th component explored in the weight-biased exploration process & Page~\pageref{def:wroot} \\
$\wst_{i}$ & time at which the exploration of $\wvx_{i}$ ends & Page~\pageref{def:wroot} \\
$\wroot_{i}$ & root vertex of $\wvx_{i}$. 
& Page~\pageref{def:wroot}\\
$(\lgclw{i},i\ge 1)$ & components of $G_n$ ordered by size & Page~\pageref{def:lgclw}\\ 
$(\hvclw{i},i\ge 1)$ & components of $G_n$ ordered by weight & Page~\pageref{def:hvclw}\\ 
$(\lgclw{i}_T,i\ge 1)$ & components of $G_n$ that are fully explored by time $Tn^{2/3}$ in weight-biased exploration, ordered by size & Page~\pageref{def:lgclwT}\\ 
$(\hvclw{i}_T,i\ge 1)$ & components of $G_n$ that are fully explored by time $Tn^{2/3}$ in weight-biased exploration, ordered by weight & Page~\pageref{def:hvclwT} \\ 
$(\swcount{t},t\ge 0)$ & the number of vertices explored in the weight-biased exploration, up to time $t$ & Page~\pageref{eq:ntilde_hat_def} \\

\hline
$\hat{G}_n$ & $\Wn$-distributed graph, constructed via the \em{size-biased exploration process} &  Page~\pageref{def:Gn hat}\\
$\svx_i$ & $i$th component explored in the size-biased exploration & Page~\pageref{def:Vk hat}\\
$\sst_{i}$ & time at which the exploration of $\svx_{i}$ ends & Page~\pageref{def:tau k hat}\\
$\sroot_{i}$ & root vertex of $\svx_{i}$ & Page~\pageref{def:sroot}\\
$(\lgcls{i},i\ge 1)$ & components of $\hat{G}_n$ ordered by size & Page~\pageref{def:lgcls}\\ 
$(\hvcls{i},i\ge 1)$ & components of $\hat{G}_n$ ordered by weight & Page~\pageref{def:hvcls}\\ 
$(\lgcls{i}_T,i\ge 1)$ & components of $\hat{G}_n$ that are fully explored by time $Tn^{2/3}$ in size-biased exploration, ordered by size & Page~\pageref{def:lgclsT} \\ 
$(\hvcls{i}_T,i\ge 1)$ & components of $\hat{G}_n$ that are fully explored by time $Tn^{2/3}$ in size-biased exploration, ordered by weight & Page~\pageref{def:hvclsT}\\ 
$(\sscount{t},t\ge 0)$ & the number of vertices explored in the size-biased exploration, up to time $t$ & Page~\pageref{eq:ntilde_hat_def} \\

\hline
$(S_n(t),t\ge 0)$ & real-valued stochastic process used to recover the components of $G_n$ and $\hat{G}_n$ & Page~\pageref{Sn} \\
$(N(t),t\ge 0)$ & the number of $i \in [n]$ with $J_i \le t$ & Page~\pageref{Nt def}\\

\hline
$\clvec$ & $\big(n^{-2/3}\big(w(\hvcls{i}),|\hvcls{i}|,w(\lgcls{i}),|\lgcls{i}|),i\in[k] \big)$ & Page~\pageref{def: W-vector} \\
$\clvec_T$ & $\big(n^{-2/3}\big(w(\hvcls{i}_T),|\hvcls{i}_T|,w(\lgcls{i}_T),|\lgcls{i}_T|),i\in[k] \big)$ & Page~\pageref{def: W-vector} \\
$\exvec_T$ & $\big( (\excsize{i}_T, \excsize{i}_T,\excsize{i}_T,\excsize{i}_T ),i\in[k]\big)$ & Page~\pageref{def: W-vector} \\
$\exvec$ & $\big( (\excsize{i}, \excsize{i},\excsize{i},\excsize{i} ),i\in[k]\big)$ & Page~\pageref{def: W-vector} \\

\hline
$(\sbm(t),t\geq 0)$ & standard Brownian motion & Page~\pageref{def:sbm} \\
$(\bm(t),t\ge 0)$ & inhomogeneous Brownian motion & Page~\pageref{defW} \\
$(\exc{i}(f),i\ge 1)$ & excursions of a function $f$, ordered by length &  Page~\pageref{def:exc} \\
$(\excsize{i}(f),i\ge 1)$ & lengths of excursions of $f$, in decreasing order &  Page~\pageref{def:fexcsize} \\
$(\excsize{i},i\ge 1)$ & lengths of excursions of $\bm$, in decreasing order &  Page~\pageref{def:excsize} \\

$(\exc{i}_T,i\ge 1)$ & excursions of $\bm|_{[0,T]}$, ordered by length & Page~\pageref{def:excT} \\

$(\excsize{i}_T,i\ge 1)$ & lengths of excursions of $\bm|_{[0,T]}$, in decreasing order &  Page~\pageref{def:excsizeT} 
\end{longtable}
\end{center}

\section{Proof of \cref{samelaw}}
\label{S:proof}
In this section we present the full proof of Theorem~\ref{samelaw}, \emph{modulo} some technical results that we defer to later sections. In Section~\ref{Sub:product topology}, we prove the joint convergence of the ranked vectors of component masses and sizes in the product topology; in Section~\ref{sec: l2} we strengthen the result to convergence in $\ell^2$.
First, however, we formally define the limit object.

For $T>0$, let $\mathbb{D}[0,T]$ denote the set of \cadlag\ functions $f:[0,T] \to \R$
equipped with the Skorohod $J_1$-topology (see, for example, \cite[Chapter 3]{billingsley1968cpm}).
Define the space $\mathbb{D}[0,\infty)$ accordingly. 

Given a function $f\in\mathbb{D}(I)$ with $I=[0,T]$ or $I=[0,\infty)$ and $x\in I$,
write 
\begin{equation}
\label{run inf}
  \underline{f}(x)\coloneqq \inf_{y\in[0, x]} f(y).
\end{equation}
\begin{dfn}[excursion]
Given a c\`adl\`ag function $f:I \to \R$, where $I \subseteq [0,\infty)$ is an interval, an {\em excursion} of $f$ is a pair $(g,d)$ with $g,d \in I$ and $0 \le g < d$, such that
    \begin{equation*}
        \min\{f(g-),f(g)\} = \underline{f}(g) = \underline{f}(d) = \min\{f(d-),f(d)\},
    \end{equation*}
    and
    \[f(x)>\underline{f}(d)\  
        \text{for all}\  x\in(g,d).
    \]
\end{dfn}

Let 
\[\Dt := \big\{f\in \mathbb{D}[0,T]:\, \forall t\in [0,T], f(t-)\le f(t)\big\}\]
be the subset of $\mathbb{D}[0,T]$ consisting of elements of $\mathbb{D}[0,T]$ with no negative jumps. 
For $f\in\Dt$ and $i\geq 1$, let 
\[
\label{def:exc}
\exc{i}(f)=(g^{(i)}(f),d^{(i)}(f))
\]
denote the $i$th longest excursion of $f$, breaking ties by increasing order of start time, where $g^{(i)}(f)$ and $d^{(i)}(f)$ are the respective start and end times of this excursion. Then define 
\[
\excsize{i}(f)  = d^{(i)}(f)-g^{(i)}(f)\, .\label{def:fexcsize}
\]
If $f$ has fewer than $i$ excursions then set $E^{(i)}(f)=\varnothing$ and $g^{(i)}(f)=d^{(i)}(f)=0$. We likewise define $E^{(i)}(f)$ for $f \in \mathbb{D}[0,\infty)$, with the convention that $E^{(i)}(f)=\varnothing$ for all $i\in\N$ in the case that $f$ has excursions of unbounded length.

 Recall that for $t\ge 0$,
\[
\bm(t) = \bm(t; \cw, \mu,\mu')=\sqrt{\frac{\mu'}{\mu}}\sbm(t) + \cw t - \frac12 \frac{\mu'}{\mu^2}t^2,
\]
for $\cw, \mu, \mu'$ as in \labelcref{eq:puvdef}-\labelcref{third moment} and $(\sbm(t),t\geq 0)$ a standard Brownian motion.  For $i\geq 1$, let $\excsize{i}\coloneqq\excsize{i}(W)$. The vector $(\excsize{i},i\ge 1)$ describes the limit appearing in \cref{samelaw}. It is known \cite[Lemma 25]{Aldous97} that almost surely $\sum_{i \ge 1} (\excsize{i})^2<\infty$, which we will use in what follows.

\subsection{Convergence in the product topology}
\label{Sub:product topology}
We begin by proving convergence in the product topology. Hereafter, we assume that $(\wtv,n\geq 1)$ is a sequence of weight vectors satisfying \labelcref{maximal weight}--\labelcref{third moment}. Let us recall from \eqref{Sn} that  
$\newrw(t) = -t + \sum_{u\in[n]} \wt{u}\I{\jump_u\leq t}$ for $t \ge 0$, and from \eqref{defW} the definition of $W$. 

\begin{prop}\label{prop:conv_S_n}
    Let $\scaledrw(t)\coloneqq n^{-1/3}S_n(tn^{2/3})$. Then it holds that 
    \[(\scaledrw(t), t\ge 0)\convdistn (W(t),t\ge 0) \]
    uniformly on compact sets. 
\end{prop}
Proposition~\ref{prop:conv_S_n} is proved using standard arguments from stochastic analysis, which we defer to \cref{sec: comment}.

For
$t\geq 0$, let
\begin{equation}
\label{eq:cchatdefs}
\wc(t) \coloneqq \sup\{i:\wst_i \leq t\}\quad\text{and}\quad \scc(t) \coloneqq \sup\{i:\sst_i \leq t\}
\end{equation}
 be the number of components explored up to time $t$ in the weight- and size-biased exploration processes, respectively. Equivalently, if we view a component as discovered at the same moment that the preceding component is fully explored, then $\wc(t)+1$ and $\scc(t)+1$ are the number of components discovered up to time $t$ in the two exploration processes.

\begin{prop}\label{prop: comps explored}
It holds that 
\begin{align*}
\left(n^{-1/3} \wc(t n^{2/3} ), t\ge 0\right)&\convdistn \left(- \inf_{0\le u\le t } W(u), t\ge 0\right) \text{ and }\\
\left(n^{-1/3} \scc(tn^{2/3}), t\ge 0\right) &\convdistn \left(- \frac{1}{\mu}\inf_{0\le u \le t } W(u), t\ge 0\right)
\end{align*}
uniformly on compact sets. Moreover, both of the previous convergences hold jointly with the one in Proposition~\ref{prop:conv_S_n}. 
\end{prop}

In reading what follows, it will be useful to recall the definitions of $J_i,R_i$ and $\hat{R}_i$ from Section~\ref{defining St}, as well as those of $(\tau_i:\,i \in [K_n])$ and $(\hat{\tau}_i:\,i \in [\hat{K}_n])$ from \eqref{eq:tauk} and \eqref{eq:taukhat}. 
 \begin{lemma}\label{roots eps dense}
    We have
    \begin{align*}
    &\Big(\frac{1}{n^{1/3}}\sum_{i=1}^{\wc(tn^{2/3})} \wt{\wroot_i}, t\ge 0\Big)\convdistn \Big(-\inf_{0\le u\le t} W(u), t\ge 0\Big) \quad\text{and}\\
    &\Big(\frac{1}{n^{1/3}}\sum_{i=1}^{\scc(tn^{2/3})} \wt{\sroot_i}, 
    t\ge 0\Big)\convdistn \Big(-\inf_{0\le u\le t} W(u), t\ge 0\Big)
    \end{align*}
    uniformly on compacts, and each holds jointly with the convergence in Proposition~\ref{prop:conv_S_n}.  
\end{lemma} 

 \begin{lemma}\label{roots dont matter}
For any $t\geq 0$, the scaled sums
\begin{align*}
&\frac{1}{n^{1/3}}\sum_{i=1}^{\wc(tn^{2/3})} \wt{\wroot_i}\I{J_{\wroot_i }\le t n^{2/3}},  
    \quad
\frac{1}{n^{1/3}}\sum_{i=1}^{\scc(tn^{2/3})}\wt{\sroot_i}\I{J_{\sroot_i} \le t n^{2/3}},\\
&\frac{1}{n^{1/3}}\sum_{i=1}^{\wc(tn^{2/3})} \I{J_{\wroot_i }\le t n^{2/3}}  
    \quad
    \mbox{and}
    \quad
\frac{1}{n^{1/3}}\sum_{i=1}^{\scc(tn^{2/3})}\I{J_{\sroot_i} \le t n^{2/3}} 
\end{align*}
all converge in probability to $0$ 
    as $n\to \infty$. 
\end{lemma}

Proposition~\ref{prop: comps explored} and Lemmas \ref{roots eps dense} and \ref{roots dont matter} are proved in \cref{root wt proofs}.

For $T>0$ and integer $i \ge 1$, write 
\[
\exc{i}_T=\exc{i}_T(\bm) = (g^{(i)}_T,d^{(i)}_T) =(g^{(i)}(\bm|_{[0,T]}),d^{(i)}(\bm|_{[0,T]}))\label{def:excT}
\]
for the $i$th longest excursion of $\bm$ that has ended by time $T$, and $\excsize{i}_T = \excsize{i}(\bm|_{[0,T]})=d^{(i)}_T-g^{(i)}_T$ for its length.\label{def:excsizeT} In what follows, we will use that $\excsize{i}_T>0$ almost surely for all $i$ and $T>0$; this is immediate from standard facts about Brownian motion. 

As a step toward proving that
\[
(n^{-2/3}w(\hvcls{i}),i\geq 1) \convdistn  (\excsize{i},i\geq 1),
\]
we first show that 
\[
(n^{-2/3}w(\hvcls{i}_T),i\geq 1) \convdistn  (\excsize{i}_T,i\geq 1).
\]

Recall 
from Section~\ref{defining St}, and in particular from \eqref{eq:taukhat},  
that the components of $\hat{G}_n$ can be extracted from $S_n$ and the root set $(\sroot_l,l\in [\hat{K}_n])$ by considering the times $(\sst_k)_{k\ge 1}$ at which $S_n$ first hits the levels $\left(-\sum_{l=1}^{k} w_{\sroot_l}\I{J_{\sroot_l} > t}\right)_{k \ge 1}$. This motivates the following definitions. 

\begin{dfn}[Cutoff levels] For $k\in \N$, a set of {\em cutoff levels} is a set of real numbers $\alpha=\{\alpha_0, \alpha_1, \dots, \alpha_k\}$ satisfying $0=\alpha_0\le \alpha_1\le\dots\le \alpha_k$. 
For $\varepsilon>0$, the set of cutoff levels $\alpha=\{\alpha_0, \alpha_1, \dots, \alpha_k\}$ is 
{\em $\veps$-dense} in a subset $I\subset\mathbb{R}$ if $I\subseteq \bigcup_{i=0}^k(\alpha_i - \veps, \alpha_i + \veps)$. 
If we refer to a set $\alpha=\{\alpha_0,\ldots,\alpha_k\}$ of cutoff levels as $\veps$-dense without specifying a set $I$, we mean that it is $\veps$-dense in $[0,\alpha_k]$.
\end{dfn}

For $T> 0$, $h\in \Dt$ such that $h(0)=0$ and $z>0$, let
\begin{equation}
  \label{tauhz}  
\tau(h,z) \coloneqq \inf\{t\geq 0:h(t)=-z\}\wedge T.
\end{equation}
For a set of cutoff levels $\alpha = \{\alpha_0,\ldots,\alpha_k\}$ and $i \in [k]$ such that $\tau(h,\alpha_i)<T$, we say that the ordered pair $(\tau(h,\alpha_{i-1}),\tau(h,\alpha_i))$ is an {\em $(h,\alpha)$-excursion.} 
For $j\ge 1$, let $g^{(j,\alpha)}(h)$ and $d^{(j,\alpha)}(h)$ be the left and right endpoints of the $j$th longest $(h,\alpha)$-excursion, breaking ties by increasing value of the left endpoint. (If the number of $(h,\alpha)$-excursions is less than $j$ then set $g^{(j,\alpha)}(h)=d^{(j,\alpha)}(h)=0$.) We write $\aexc{j}(h)\coloneqq (g^{(j,\alpha)}(h),d^{(j,\alpha)}(h))$ and $\excsize{j,\alpha}(h)=d^{(j,\alpha)}(h)-g^{(j,\alpha)}(h)$. 

The following lemma states that, for a sequence of processes that converges in distribution to $(\bm(t), 0\le t \le T)$, the ordered excursions demarcated by  sufficiently dense sets of cutoff levels closely approximate the ordered excursions of $(\bm(t),0 \le t \le T)$ as $n\to \infty$.

\begin{lemma}\label{lem:eps_deldisguise}
Let $T>0$. 
For $n\geq 1$, let $\overline{X}_n$ be a stochastic process with paths in $\mathbb{D}^+[0,T]$ that converges uniformly to
$(\bm(t),t \in [0,T])$ in distribution as $n\to \infty$. For $n \geq 1$, let $\alpha^{(n)}\coloneqq \{\alpha^{(n)}_0,  \dots, \alpha^{(n)}_{m_n}\}$ and $\beta^{(n)}\coloneqq \{\beta^{(n)}_0,  \dots,\beta^{(n)}_{m'_n}\}$ be two sets of cutoff levels such that for any $\veps>0$, 
\begin{equation}\label{within eps of inf}
\liminf_{n\to\infty}\Pp\big( - \min\{\alpha_{m_n}^{(n)}, \beta_{m'_n}^{(n)}\} < \inf \overline{X}_n+ \veps\big) > 1- \veps ,
\end{equation}
and for all $n$ large enough, $\alpha^{(n)}$ and $\beta^{(n)}$ are $\veps$-dense. Then  
\[
\Big(\big(E^{(i,\alpha^{(n)})}(\overline{X}_n) ,  E^{(i,\beta^{(n)})}(\overline{X}_n)\big), i \geq 1\Big)
\convdistn \Big(\big(\exc{i}_T(\bm),\exc{i}_T(\bm)\big),i \ge 1\Big)
\]
jointly in the product topology.
\end{lemma}
We prove \cref{lem:eps_deldisguise} in \cref{proof of disguised delta eps}, by appealing to a deterministic result that holds for convergent sequences with a limit in a suitable class of ``good'' functions and using the Skorohod representation theorem.

For $T>0$, we let $\hvclw{i}_T$\label{def:hvclwT} and $\lgclw{i}_T$\label{def:lgclwT} be the vertex sets of the $i$th heaviest and $i$th largest components, respectively, whose explorations in the weight-biased process have been completed by time $Tn^{2/3}$. Similarly, we let $\hvcls{i}_T$\label{def:hvclsT} and $\lgcls{i}_T$\label{def:lgclsT} be the vertex sets of the $i$th heaviest and $i$th largest components, respectively, whose explorations in the size-biased process have been completed by time $Tn^{2/3}$. Since the processes explore components in weight-biased and in size-biased order, respectively, it need not hold that $\hvclw{i}_T \eqdist \hvcls{i}_T$ or that $\lgclw{i}_T\eqdist \lgcls{i}_T$. However, the next proposition implies that, for fixed $i$, their weights are with high probability close to one another when $n$ is large.
We denote by $g^{(n, i)}_T$ and $d^{(n, i)}_T$ the respective starting and ending times of the exploration of $\hvclw{i}_T$ in the weight-biased exploration. We likewise define $\hat g^{(n, i)}_T$ and $\hat d^{(n, i)}_T$ for $\hvcls{i}_T$. 

\begin{prop}\label{joint convergence to exc kT}
It holds that 
\[
\Big(n^{-2/3}\Big(\big(\hat g^{(n, i)}_T, \hat d^{(n, i)}_T\big), \big(g^{(n, i)}_T, d^{\,(n, i)}_T\big)\Big), i\ge 1 \Big) \convdistn  \big( \big(\exc{i}_T,\exc{i}_T\big), i\ge 1 \big)
\]
in the product topology. 
\end{prop}
\begin{proof}
Recall that $\scaledrw(t)=n^{-1/3}S_n(n^{2/3}t)$.
We will prove the proposition by applying \cref{lem:eps_deldisguise} to the process $\scaledrw$ and the cutoff level sets used in the size-biased and the weight-biased exploration processes. For $n\geq 1$ and $k\leq m_n\coloneqq \scc(Tn^{2/3})$, we define 
\[
\alpha^{(n)}_k = \frac{1}{n^{1/3}}\sum_{l\in[k]} w_{\sroot_l}\I{J_{\sroot_l} > \sst_{k}}\]
and consider the set of cutoff levels $\alpha^{(n)}\coloneqq \{\alpha^{(n)}_0,  \cdots, \alpha^{(n)}_{m_n}\}$. Similarly, for $n\ge 1$ and $k\le m'_n\coloneqq \wc(Tn^{2/3})$, we define
\[
\beta^{(n)}_k = \frac{1}{n^{1/3}}\sum_{l\in[k]} w_{\wroot_l}\I{J_{\wroot_l} > \wst_{k}}
\]
and let $\beta^{(n)}\coloneqq \{\beta^{(n)}_0, \cdots, \beta^{(n)}_{m'_n}\}$.
By \eqref{maximal weight}, for any $\veps > 0$, for all $n$ large enough, $\alpha^{(n)}$ and $\beta^{(n)}$ are $\veps$-dense. 

Recall that $\sst_i$ is the time at which the exploration of $(i+1)$st component, $\svx_{i+1}$, begins in the size-biased exploration process. 
Then, if $\hvcls{k}_T=\svx_{i+1}$, we have $\hat g^{(n, k)}_T=\sst_i$ and $\hat d^{(n, k)}_T=\sst_{i+1}$.
Moreover, $\sst_i$ can be recovered from the random walk $S_n$ as
\[
\sst_{i}=\inf\Big\{t\geq 0:\newrw(t) =-\sum_{j\in [i]} \wt{\sroot_j}\I{J_{\sroot_j} > t}\Big\} = n^{2/3}\inf\left\{t\geq 0:\scaledrw(t) =-\alpha^{(n)}_i\right\}.
\] 
Recalling that $m_n=\hat c(Tn^{2/3})$, we deduce from \cref{roots eps dense} and \cref{roots dont matter} that
\begin{align*}
\alpha_{m_n}^{(n)}
& =\frac{1}{n^{1/3}}\sum_{i=1}^{m_n}\wt{\sroot_i}\I{\{J_{\sroot_i}> t\}} \\
& = \frac{1}{n^{1/3}}\sum_{i=1}^{\hat c(Tn^{2/3})}\wt{\sroot_i}-\frac{1}{n^{1/3}}\sum_{i=1}^{\hat c(Tn^{2/3})}\wt{\sroot_i}\I{\{J_{\sroot_i}\le t\}}\\
& \convdist -\inf_{t\in [0,T]}W(t)
\end{align*}
as $n\to \infty$, 
jointly with the convergence of $(\scaledrw(t))_{t\ge 0}$ to $W$. This implies that
\[
\Big(-\alpha_{m_n}^{(n)}, \inf_{t\in [0, T]}\scaledrw(t) \Big) \convdistn \Big(\inf_{t\in [0,T]}W(t), \inf_{t\in [0,T]}W(t)\Big).
\]
So we get that for $n$ large enough,
\[
\Pp\left(-\alpha_{m_n}^{(n)}<\inf_{t\in[0,T]} \scaledrw(t) +\veps\right) > 1-\veps.
\]
A similar argument yields
\[
\Pp\left(-\beta_{m'_n}^{(n)}<\inf_{t\in[0,T]} \scaledrw(t) +\veps\right) > 1-\veps.
\]
Thus, \eqref{within eps of inf} is satisfied by the process $\scaledrw$ with cutoff levels $\alpha^{(n)},\beta^{(n)}$. 
An application of~\cref{lem:eps_deldisguise} then concludes the proof.
\end{proof}
To deduce that the {\em sizes} of the heaviest components also converge under rescaling, we will study the counting process that keeps track of the number of vertices discovered up to a given time. To that end, define
\begin{equation*}\label{Nt def}
\scount{t} = \scount{t}(\wtv)\coloneqq\sum_{i\in[n]}\I{\jump_i\leq t}.
\end{equation*}

\begin{lemma}\label{LLN lemma}
We have $(n^{-2/3}\scount{n^{2/3}t}, t\ge 0)\to (t, t\ge 0)$ in probability uniformly on compacts as $n\to \infty$.  
\end{lemma}

We prove \cref{LLN lemma} in \cref{sec: comment}.

For $t\ge 0$, let us define
\begin{equation}\label{eq:ntilde_hat_def}
\swcount{t}=\scount{t}
+\sum_{i= 1}^{\wc(t)}\I{J_{\wroot_i}> t}\quad\text{and}\quad \sscount{t}=\scount{t}
+\sum_{i= 1}^{\scc(t)}\I{J_{\sroot_i}> t}, 
\end{equation}
which count precisely the number of vertices discovered in the weight- and size-biased exploration processes, respectively, up to time $t$.
In particular, since the component $\wvx_k$ (resp.\ $\svx_k$) is explored during the interval $[\wst_{k-1}, \wst_{k})$ (resp.\ $[\sst_{k-1}, \sst_{k})$), it holds that
\begin{equation}
\label{explore-times}
|\wvx_k| = \swcount{\wst_{k}-}- \swcount{\wst_{k-1}-}\quad\text{and}\quad |\svx_k| = \sscount{\sst_{k}-}- \sscount{\sst_{k-1}-},
\end{equation}
where we take $N(0-)=\hat{N}(0-)=0$ by convention.

By \cref{roots dont matter} and~Proposition~\ref{prop: comps explored}, it holds that 
\[
\swcount{tn^{2/3}} = \scount{tn^{2/3}} + O_\Pp(n^{1/3})\quad\text{and}\quad\sscount{tn^{2/3}} = \scount{tn^{2/3}} + O_\Pp(n^{1/3}),
\]
so we have the following corollary of \cref{LLN lemma}.

\begin{cor}\label{LLN 2}
It holds that for all $\veps>0,T>0$ and for $n$ large enough,
\begin{align}
\Pp\left(\sup_{t\le Tn^{2/3}} n^{-2/3}|\swcount{t} - t|>\veps\right)&<\veps, \label{Nwunifconv}\\
\Pp\left(\sup_{t\le Tn^{2/3}} n^{-2/3}|\sscount{t} - t|>\veps\right)&<\veps. \label{Ntunifconv}
\end{align}
\end{cor}

\begin{prop}\label{sizes conv T}
For all $T>0$, 
\begin{equation}\label{joint conv weight biased}
\left(n^{-2/3}\left(|\hvclw{i}_T|,w(\hvclw{i}_T),|\hvcls{i}_T|,w(\hvcls{i}_T)\right), i\ge 1\right) \convdistn  \big(\excsize{i}_T,\excsize{i}_T,\excsize{i}_T,\excsize{i}_T, i\ge 1\big)
\end{equation}
in the product topology.
\end{prop}
\begin{proof}
We will show that for any $i\geq 1$, 
\begin{equation}\label{eq:jointcon_toprove}
\left(\frac{|\hvcls{i}_T|}{w(\hvcls{i}_T)},\frac{|\hvclw{i}_T|}{w(\hvclw{i}_T)}\right)  \convprobn (1,1).
\end{equation}
Recalling that $w(\hvcls{i}_T)=\hat d^{(n, i)}_T-\hat g^{(n, i)}_T$ and $w(\hvclw{i}_T)=d^{(n, i)}_T- g^{(n, i)}_T$,
the result then follows from \cref{joint convergence to exc kT} and Slutsky's theorem.

Fix $i\geq 1$ and $\veps' > 0$. By \cref{joint convergence to exc kT} and since $\excsize{i}_T>0$ almost surely, 
we can find $\delta\in (0, 1)$ such that for all $n$ large enough, \[
\Pp\left(n^{-2/3}w(\hvcls{i}_T) < \delta \right) < \frac{\veps'}{2}.
\]
Recall that $\hat g^{(n,i)}_T$ and $\hat d^{(n, i)}_T$ are the
starting and ending times, respectively, of the exploration of $\hvcls{i}_T$, so $(\hat g^{(n,i)}_T,\hat d^{(n, i)}_T)=(\hat{\tau}_\ell,\hat{\tau}_{\ell+1})$ for some $\ell \ge 0$ with $\hat{\tau}_{\ell+1} \le Tn^{2/3}$, so
\[
|\hvcls{i}_T| = \sscount{\hat d^{(n, i)}_T-} -\sscount{\hat g^{(n, i)}_T-} \text{ and } 
w(\hvcls{i}_T) = \hat d^{(n, i)}_T-\hat g^{(n, i)}_T.
\]
Writing $\veps = \frac{\delta}{2} \cdot \veps'$, then it follows from \labelcref{Ntunifconv} that for all $n$ large enough,
\begin{align*}
&\frac{\veps'}{2}>\Pp\left(n^{-2/3}\left|(\sscount{\hat d^{(n, i)}_T}-\sscount{\hat g^{(n, i)}_T-}- (\hat d^{(n, i)}_T-\hat g^{(n, i)}_T)\right|\ge \veps \right)\\
&=\Pp\left(\left||\hvcls{i}_T|- w(\hvcls{i}_T)\right|\ge \veps n^{2/3} \right)\\
&= \Pp\left(\left| \frac{|\hvcls{i}_T|}{w(\hvcls{i}_T)} - 1\right| \ge  \frac{ \veps }{n^{-2/3} w(\hvcls{i}_T)}\right). 
\end{align*}
Then by our choice of $\delta$, it holds that for all $n$ large enough,
\[
\left| \frac{|\hvcls{i}_T|}{w(\hvcls{i}_T)} - 1\right|<\frac{\veps}{\delta}  < \veps'
\]
with probability at least $1-\veps'$. Analogously, we can show that
\[
\Pp\left(\left|\frac{|\hvclw{i}_T|}{w(\hvclw{i}_T)}-1\right|>\veps'\right)<\veps',
\]
thus obtaining the joint convergence \eqref{eq:jointcon_toprove}.
\end{proof}

Let us recall that $\exc{i}(W)=(g^{(i)},d^{(i)})$ (resp.~$\exc{i}_T(W)=(g^{(i)}_T,d^{(i)}_T)$) stands for the $i$th longest excursion of $W$ (resp.~ of $(W(t), 0\le t\le T)$), and that $\excsize{i}=d^{(i)}-g^{(i)}$,  $\excsize{i}_T=d^{(i)}_T-g^{(i)}_T$.
Note that for all $T$ large enough (specifically, for $T \ge \max\{d^{(j)}, j\in [i]\}$), we have $\excsize{i}_T=\excsize{i}$ for each $i\ge 1$. This immediately yields the following statement.
\begin{lemma}\label{brownian excs}
For all $k \ge 1$, the following holds almost surely: for all $T$ large enough, 
\[
(\excsize{1}_T,\ldots,\excsize{k}_T) =(\excsize{1},\ldots,\excsize{k}).
\]
\end{lemma}

For the next lemma, recall that $\lgcls{j}_T$ stands for the $j$th largest connected component explored up to time $Tn^{2/3}$ in the size-biased exploration.
\begin{lemma}\label{hvT equal lgT}
For $k\ge 1$, it holds that 
\[
\lim_{T\to \infty}\liminf_{n\to\infty}\Pp
\left(\hvcls{j}_T  = \lgcls{j}_T\text{ for all } j\in [k]\right) =1.
\]
\end{lemma}

\begin{proof}
Fix $k\geq 1$ and $\veps \in (0,\frac{1}{2})$. 
Since the sequence $(\excsize{i},i\ge 1)$ is strictly decreasing and converges to $0$ almost surely, we may pick $\delta>0$ and $J\in\N$ such that $\Pp\left(\excsize{k} > \delta > \excsize{J} \right) > 1-\frac{\veps}{16}$. 
By \cref{brownian excs}, we may take $T$ large enough that 
$\Pp\left(\excsize{i}_T=\excsize{i}\text{ for all } i\in [J]\right) > 1-\frac{\veps}{16}$.
Then, by \cref{sizes conv T} and since  $(\excsize{i},i\in [J])$ is strictly decreasing almost surely, 
for $n$ large enough, with  probability at least $1-\frac{\veps}{4}$, both
\begin{equation}\label{strictly decreasing}
|\hvcls{1}_T| > \cdots > |\hvcls{J}_T|
\end{equation}
and 
\begin{equation}\label{size and weight bound}
 n^{-2/3}|\hvcls{1}_T| > \delta > n^{-2/3}w(\hvcls{J}_T) 
\end{equation}
hold. 

Let $\jn$ be such that $\lgcls{1}_T = \hvcls{\jn}_T$.
We will show that for all $n$ large enough, $\Pp(\jn = 1) > 1 - {\veps}$. By \labelcref{strictly decreasing}, for all $n$ large enough, 
\[
\Pp\left(\lgcls{1}_T \ne \hvcls{1}_T\right) =\Pp\left(\jn>1\right)\le \frac{\veps}{4} + \Pp\left(\jn >J\right).
\]
Now we show that $\Pp(\jn>J)\le\frac{\veps}{2}$. 
Let $\badcomponent$ be the event that $|\lgcl{1}_T| > \delta n^{2/3} > w(\lgcls{1}_T)$. 
Then 
\begin{align}
    \Pp(\jn > J)&\le \Pp\left((\badcomponent)^c \cap \{\jn>J\}\right)+\Pp(\badcomponent) \nonumber\\
    &\le \Pp\left(|\lgcls{1}_T| \le  \delta n^{2/3} \ \text{or}\ \delta n^{2/3} \le w(\lgcls{1}_T)\,\middle|\, \jn>J \right)+\Pp(\badcomponent) \nonumber\\
    &\le \Pp\left(|\hvcls{1}_T| \le  \delta n^{2/3} \ \text{or}\ \delta n^{2/3} \le w(\hvcls{J}_T)\right)+\Pp(\badcomponent), \label{hvTequalslgT:ineq1}
\end{align}
because  $|\lgcls{1}_T|\ge |\hvcls{1}_T|$ and, if $\jn \ge J$, then $w(\hvcls{J}_T)> w(\lgcls{1}_T)$.  
The first term is at most $ \frac{\veps}{4}$ by \eqref{size and weight bound}. 
To bound the second term, recall that in the size-biased exploration process, the component $\lgcls{1}_T$ is explored during an interval of length exactly $w(\lgcls{1}_T)$, and $|\lgcls{1}_T|$ corresponds to the number of jumps of $\hat N$ that occur within this interval. 
Hence, $\Pp(\badcomponent)$ is less than the probability that 
there exists $t\in  [0,(T-\delta)n^{2/3}]$ such that $\sscount{t+\delta n^{2/3}} - \sscount{t} \ge 2\delta n^{2/3}$. Then,
\begin{align}
\Pp(\badcomponent) &\le  \Pp\left(\sup_{t\in [0,(T-\delta)n^{2/3 }]} |\sscount{t} - \sscount{t+\delta n^{2/3}}|>2\delta n^{2/3} \right) \nonumber\\
&\le \Pp\left(\sup_{t\in [0,Tn^{2/3}]} n^{-2/3}|\sscount{t} - t|>\min\Big\{\frac{\veps}{4},\frac{\delta}{2}\Big\} \right) 
 \le \frac{\veps}{4}\, , \label{hvTequalslgT:ineq2}
\end{align}
the last inequality holding for all $n$ sufficiently large by \labelcref{Ntunifconv}.
Combining \eqref{hvTequalslgT:ineq1} and \eqref{hvTequalslgT:ineq2} yields that
$\Pp(\jn > J) < \frac{\veps}{2}$, so that 
\[\Pp\Big(\lgcls{1}_T = \hvcls{1}_T\Big)>1-\veps.\]

By a similar argument, we find that for all $n$ large enough, 
\begin{equation*}
\Pp\left(\hvcls{j}_T  = \lgcls{j}_T\text{ for all $j\in [k]$}\right) >  1- k \veps.
\end{equation*}
Since $\veps$ can be taken arbitrarily small, this completes the proof.
\end{proof}

For $k\ge 1$ and $T>0$, we define the vectors 
\begin{align}
\notag
&\clvec_T = \big(n^{-2/3}\big(w(\hvcls{i}_T),|\hvcls{i}_T|,w(\lgcls{i}_T),|\lgcls{i}_T|),i\in [k] \big),\\ \notag
&\clvec = \big(n^{-2/3}\big(w(\hvcls{i}),|\hvcls{i}|,w(\lgcls{i}),|\lgcls{i}|),i\in [k] \big),\\ \notag
&\exvec_T=\big( (\excsize{i}_T, \excsize{i}_T,\excsize{i}_T,\excsize{i}_T ),i\in [k]\big)\quad\text{ and }\\ \label{def: W-vector}
&\exvec=\big( (\excsize{i}, \excsize{i},\excsize{i},\excsize{i} ),i\in [k]\big).
\end{align}

\begin{lemma}\label{conv of components}
For $k\ge 1$ it holds that $\clvec_T \convdist \exvec_T$ as $n\to \infty$.
\end{lemma}
\begin{proof}
Fix $\veps \in(0,1), k\ge 1$, and a continuous function $f:\mathbb{R}^{4k}\to [0, 1]$.
For $T>0$ and $n\ge 1$, define the event $\cA_T^{(n,k)} = \big\{\hvcls{j}_T  = \lgcls{j}_T\text{ for all $j\in [k]$} \}$. By \cref{hvT equal lgT}, we may choose $T>0$ such that for all $n$ large enough, 
\begin{equation}\label{take n large}
\Pp(\cA_T^{(n,k)}) > 1-\frac{\veps}{4}.
\end{equation}

By \cref{sizes conv T}, it holds that for all $n$ large enough, 
\begin{equation}\label{conv f}
\left | \mathbb{E}\left[f \left(n^{-2/3}(w(\hvcls{i}_T),|\hvcls{i}_T|,w(\hvcls{i}_T),|\hvcls{i}_T|),i\in [k] \right) - f\left(\exvec_T \right) \right] \right| < \frac{\veps}{2}.
\end{equation}
Then for all $n$ large enough,
\begin{align*}
\left|\mathbb{E}\left[f\left(\clvec_T\right) - f \left(\exvec_T\right)\right]\right| &\le \left|\mathbb{E}\left[ \left(f\left(\clvec_T \right) - f \left(\exvec_T\right)\right)\cdot\I{\cA_T^{(n,k)}}\right]\right|\nonumber\\
& + \left|\mathbb{E}\left[ \left(f\left(\clvec_T \right) - f \left(\exvec_T\right)\right)\cdot\I{(\cA_T^{(n,k)})^c}\right]\right|\nonumber\\
    &\le \frac{\veps}{2} + 2\cdot\Pp\Big((\mathcal{A}_T^{(n,k)})^c\Big) < \veps \quad\text{by \eqref{take n large} and \eqref{conv f}.}
\end{align*}
Hence, for any continuous function $f:\R^{4k}\to [0, 1]$, it holds that
\[
\mathbb{E}\left[f \left(\clvec_T \right) \right] \to \mathbb{E}\left[f \left(\exvec_T \right) \right]
\]
as $n\to\infty$, so $\clvec_T \to \exvec_T$ in distribution as $n\to \infty$.
\end{proof}
The next lemma is one of the crucial results of the section; it states that both the largest and heaviest components of $\hat{G}_n$ are very likely to have been discovered by time $Tn^{2/3}$, provided $T$ is sufficiently large.
\begin{lemma}\label{rhino}
For any $k\ge 1$, 
\begin{align}
\lim_{T\to \infty}\liminf_{n\to \infty}\Pp\left( (\lgcls{i}_T,i\in [k]) = (\lgcls{i},i\in[k])\right)&=1\mbox{ and} \label{rhino size}\\
\lim_{T\to \infty}\liminf_{n\to \infty} \Pp\left((\hvclw{i}_T,i\in[k]) = (\hvclw{i},i\in[k])\right)&=1.\label{rhino weight}
\end{align}
\end{lemma} 

\begin{proof}

Fix $\veps\in (0,1)$ and $k\ge 1$. By \cref{brownian excs}, there exists $T$ large enough such that for all $n$ large enough, we have
\begin{equation} \label{k largest brownian excs}
    \Pp\big(\excsize{i}_{T} = \excsize{i},i\in [k]\big) \ge 1-\frac{\veps}{4}.
\end{equation}
Since $\excsize{k} > 0$ almost surely, 
there exists $\delta>0$ such that $\Pp(\excsize{k}<\delta)<\frac{\veps}{4}$, 
which along with \eqref{k largest brownian excs} yields that
\begin{equation}\label{excursions lower bound}
\Pp\left(\excsize{k}_T<\delta \right)<\frac{\veps}{2}.
\end{equation}
From \cref{conv of components}, it follows that for $n$ large enough, 
\[
\Pp\left(|\lgcls{k}_T|<\delta n^{2/3}\right) < \veps.
\] 
Hence, to prove \eqref{rhino size}, it suffices to prove that for $n$ and $T$ large enough, with probability greater than $1-\veps$, all components of size at least $\delta n^{2/3}$ have been fully explored by time $Tn^{2/3}$ by the size-biased exploration process.

Let $\badr$ be the event that there is a component with size at least $\delta n^{2/3}$ whose exploration has not been completed by time $Tn^{2/3}$ by the size-biased exploration process. It will be convenient for our proofs to consider the time that a component is discovered rather than fully explored. This motivates the following definitions: let $\badrt$ be the event that there is a component with size at least $\delta n^{2/3}$ that has not been discovered by time $\tfrac{1}{2}Tn^{2/3}$ by the size-biased exploration process, and let $\badrh$ be the event that there is a component that is discovered by time $\tfrac 12 Tn^{2/3}$ but not fully explored by time $Tn^{2/3}$ (in which case $\hat{c}(\tfrac 12 Tn^{2/3})=\hat{c}(Tn^{2/3})$).
Then
it holds that 
\[\Pp(\badr)\le \Pp(\badrt)+\Pp(\badrh).\]
We will show that for $T$ large enough and $n$ large enough, both $\Pp(\badrt)$ and $\Pp(\badrh)$ are at most $\frac{\veps}{2}$. 

By Claim~\ref{claim:inf_W}, if $T$ is sufficiently large, then 
\[
\Pp\left(\inf_{0\le s \le T/2} W(s)-\inf_{0\le s \le 2T}W(s)<1\right)<\frac{\veps}{4}.
\]
Also, by Proposition~\ref{prop: comps explored},
$(n^{-1/3}\scc(tn^{2/3}),t \ge 0) \convdist (-\frac{1}{\mu}\inf_{0\le s \le t} W(s),t \ge 0)$ as $n \to \infty$. So for $T$ large enough, for all $n$ sufficiently large,
\[
\Pp(\badrh)=\Pp\big(\scc(\tfrac 12 Tn^{2/3})=\scc(T n^{2/3})\big)<\frac{\veps}{2}.
\] 
We  bound $\Pp(\badrt)$ using the following idea. If there is a large component that has not been discovered by time $\tfrac 12Tn^{2/3}$ by the size-biased exploration process, then it is likely to be discovered between time $\tfrac 12Tn^{2/3}$ and $2Tn^{2/3}$. However, using the convergence of $\scaledrw(t)$ to $\bm(t)$, we will show that we are unlikely to find any large components between time $\tfrac 12 Tn^{2/3}$ and $2Tn^{2/3}$.

Let $\sawr$ be the event that a component of size at least $\delta n^{2/3}$ is discovered between time $\tfrac 12 Tn^{2/3}$ and $2Tn^{2/3}$ of the size-biased exploration process. Also, fix $K>0$ to be chosen later in the proof, and define the events 
\begin{align*}
    \sawAr &\!=\!\{\text{the first $Kn^{1/3}$ components discovered after time $\tfrac 1 2 Tn^{2/3}$ are smaller than $\delta n^{2/3}$} \},\\
    \sawBr &\!=\!\{\text{fewer than $Kn^{1/3}$ components are discovered between times $\tfrac 12 Tn^{2/3}$ and $2Tn^{2/3}$}\} \\
    & \!=\!
    \{
    \scc(2Tn^{2/3})-\scc(\tfrac 12 Tn^{2/3})< Kn^{1/3}
    \}\, ,
\end{align*}
where we recall from \eqref{eq:cchatdefs} that $\scc(t)+1$ is the number of components discovered by time $t$ in the size-biased exploration process. It is immediate that $(\sawr)^c\subseteq \sawAr \cup \sawBr$, so 
\begin{align}
\Pp(\badrt)&\le \Pp(\sawAr\cap \badrt)+\Pp(\sawBr\cap \badrt)+\Pp(\sawr\cap  \badrt)\nonumber\\
&\le \Pp(\sawAr\mid \badrt)+\Pp(\sawBr)+\Pp(\sawr).\label{eq:bnt_bd}
\end{align}
So to prove that $\Pp(\badrt)\le \frac{\veps}{2}$, it suffices to show that for $T$ large enough and $K$ chosen appropriately, all three terms in the right-hand side are at most $\frac{\veps}{6}$. We bound the terms in their order of appearance. 

Each time that we begin exploring a new component, we choose each remaining component with probability proportional to its size. If there is a component with size at least $\delta n^{2/3}$ that is undiscovered by time $\tfrac 12 Tn^{2/3}$, then each time when we begin exploring a new component between times $\tfrac 12 Tn^{2/3}$ and $2Tn^{2/3}$, the probability for that particular component to be chosen is at least 
\[
\frac{\delta n^{2/3}}{n-\sscount{\tfrac 12 Tn^{2/3}}} \ge \frac{\delta n^{2/3}}{n}=\delta n^{-1/3},
\]
where we recall from \eqref{eq:ntilde_hat_def} that $\hat{N}(t)$ is the number of vertices seen in the size-biased exploration process up to time $t$.

Thus, using that $\log(1-x)\le -x$ for all $x\ge 0$, we find that
\[
\Pp(\sawAr\mid\badrt) \leq \left( 1-\delta n^{-1/3} \right)^{Kn^{1/3}}\le e^{-\delta K}.
\]
We may thus pick $K=K(\delta)$ large enough so that $\Pp(\sawAr|\badrt) < \frac{\veps}{6}$. We emphasize that this bound holds for all $n$ large enough that $1-\delta n^{-1/3}>0$ and for all $T>0$.

Next, by Claim~\ref{claim:inf_W} for all $T$ large enough, it holds that 
\[\Pp\left(\inf_{0\le s \le T/2} W(s)-\inf_{0\le s \le 2T}W(s)<K\right)<\frac{\veps}{12} .\]
Using the convergence in distribution from Proposition~\ref{prop: comps explored} as above, it follows that for $n$ large enough, 
\[\Pp(\sawBr)=\Pp\left(\scc(2Tn^{2/3})-\scc(\tfrac 12 Tn^{2/3})<Kn^{1/3}\right)<\frac{\veps}{6} .\]
We now bound $\Pp(\sawr)$. We will first show that, with probability at least $1-\frac{\veps}{12}$, for $n$ large enough, any component discovered by time $2Tn^{1/3}$ is fully explored by time $3Tn^{1/3}$. Indeed, again, by the convergence in distribution of $n^{-1/3}\scc(tn^{2/3})$   to $-\inf_{0\le s \le t} W(s)$, we may pick $T$ large enough and $n$ large enough that $\Pp(\scc(3Tn^{2/3})>\scc(2T n^{2/3}))>1-\frac{\veps}{12}$.

Let $\excsize{1}_{(a,b)}$ denote the length of the largest excursion of $\bm$ that has both its endpoints in the interval $(a,b)$. 
We similarly define $\lgcls{1}_{(a,b)}$ to be the largest component discovered after time $an^{2/3}$ and fully explored by time $bn^{2/3}$ by the size-biased process. From \cref{brownian excs} and the fact that the ordered excursions of $\bm$ are square-summable (\cite[Lemma 25]{Aldous97}, and noted just before Section~\ref{Sub:product topology}), we can choose $T$ large enough such that 
\begin{equation} \label{W large exc in T to 2T}
\Pp\left(\excsize{1}_{(\frac{1}{2}T,3T)} > \delta\right) < \frac{\veps}{24}. 
\end{equation}
Then, by the same arguments used in the proof of Proposition~\ref{sizes conv T} and Lemma~\ref{hvT equal lgT}, we can deduce from \cref{joint convergence to exc kT} and \cref{LLN 2} that
\[n^{-2/3}\big|\lgcls{1}_{(\frac{1}{2}T,3T)}\big|\convdistn\excsize{1}_{(\frac{1}{2}T,3T)}\,,\]
so we may pick $n$ large enough that
\[\Pp\left(\big|\lgcls{1}_{(\frac{1}{2}T,3T)}\big| > \delta n^{2/3}\right) < \frac{\veps}{12}. \]
Thus, 
\[\Pp(\sawr)\le\Pp\big(\scc(3Tn^{2/3})=\scc(2T n^{2/3})\big)+\Pp\left(\big|\lgcls{1}_{(\frac12 T,3T)}\big| > \delta n^{2/3}\right)<\frac{\veps}{6}. \]
and using the three bounds in \eqref{eq:bnt_bd} we deduce that $\Pp(\badrt)<\frac{\veps}{2}$ for $n$ sufficiently large. 
This concludes the proof for \eqref{rhino size}.

The same argument with some slight modifications will prove \eqref{rhino weight}. We introduce the respective analogues of $\sawAr$ and $\badrt$: let $\mathcal A^{(n,\mathrm{w})}_T$ be the event that in the weight-biased exploration the first $Kn^{1/3}$ components discovered after time $\frac12 Tn^{2/3}$ have weights smaller than $\delta n^{2/3}$, and let $\tilde{\mathcal B}^{(n,\mathrm{w})}_T$ be the event that there is a component with weight at least $\delta n^{2/3}$ that has not been discovered by time $\frac12Tn^{2/3}$ in the same exploration process.

We only argue in the case that $n$ is large enough. Given $\tilde{\mathcal B}^{(n,\mathrm{w})}_T$, each time we choose a new component for exploration after time $\tfrac 1 2 T n^{2/3}$, the probability of choosing a component of weight no smaller than $\delta n^{2/3}$ is at least
\[ 
\frac{\delta n^{2/3}}{\sum_{i\in[n]} \wt{i}}\ge  \frac{\delta n^{2/3}}{2\mi n} = \frac{\delta}{2\mi n^{1/3}}, 
 \]
where the inequality follows from  \eqref{first moment}, since for $n$ large enough, $\sum_{i\in [n]} \wt{i} \le 2\mi n$. 
It follows by the same logic as above that $\Pp(\mathcal A^{(n,\mathrm{w})}_T\mid\tilde{\mathcal B}^{(n,\mathrm{w})}_T)$ is bounded from above by
\[
 \left(1-\frac{\delta}{2\mi n^{1/3}}\right)^{Kn^{1/3}}\le e^{-K\delta/2\mi}<\frac{\veps}{6},
\]
where the final inequality holds for large enough $K$. Since the rest of the proof is essentially identical to that of \eqref{rhino size}, we omit the details.
\end{proof}

\begin{lemma}\label{new aldous}
It holds that
\begin{equation}\label{conv of weights w}
\big( n^{-2/3} (w(\hvclw{i}),|\hvclw{i}| ) ,i\ge 1 \big)\convdistn \big((\excsize{i},\excsize{i}),i\ge 1 \big)
\end{equation}
and 
\begin{equation}\label{conv of weights s}
\big( n^{-2/3} (w(\hvcls{i}),|\hvcls{i}| ) ,i\ge 1 \big)\convdistn \big((\excsize{i},\excsize{i}),i\ge 1 \big)
\end{equation}
in the product topology.
\end{lemma}
\begin{proof}

Since $\big((w(\hvcls{i}),|\hvcls{i}|),i\ge 1\big)\eqdist \big((w(\hvclw{i}),|\hvclw{i}|),i\ge 1\big)$, it suffices to prove \labelcref{conv of weights w}.

Fix $k \ge 1$, and for $T>0$ define the vectors 
\begin{align*}
    \mathrm{U}^{(n,k)}_T &= \big(n^{-2/3}(w(\hvclw{i}_T),|\hvclw{i}_T|),i\in [k]\big),\\
    \mathrm{U}^{(n,k)} &= \big(n^{-2/3}(w(\hvclw{i}),|\hvclw{i}|),i\in [k]\big),\\
    \mathrm{U}^{(k)}_T &= \big((\excsize{i}_T,\excsize{i}_T),i\in [k]\big),\\
    \mathrm{U}^{(k)} &= \big((\excsize{i},\excsize{i}),i\in [k]\big).
\end{align*}
We use the principle of accompanying laws (see, for instance \cite[Theorem 9.1.13]{Stroock_2024}) to prove that $\mathrm{U}^{(n,k)}\convdist \mathrm{U}^{(k)}$ as $n\to \infty$. By \cref{rhino}, we have
\[
\lim_{T\to \infty}\lim_{n\to \infty} \Pp\left((\hvclw{i}_T,i\in [k]) = (\hvclw{i},i\in [k])\right) = 1,\]
which implies that 
\[
\lim_{T\to \infty}\limsup_{n\to \infty}\Pp\left(\mathrm{U}^{(n,k)}_T\ne \mathrm{U}^{(n,k)}\right)=0.
\]
By Proposition~\ref{sizes conv T}, $\mathrm{U}^{(n,k)}_T\convdist \mathrm{U}^{(k)}_T$ as $n\to \infty$.
Thus, by the principle of accompanying laws, there exists a random variable $\tilde{\mathrm{U}}^{(k)}$ such that $\mathrm{U}^{(k)}_T\convdist \tilde{\mathrm{U}}^{(k)}$  as $T\to \infty$ and $\mathrm{U}^{(n,k)}\convdist \tilde{\mathrm{U}}^{(k)}$ as $n\to \infty$. But by \cref{brownian excs}, $\mathrm{U}^{(k)}_T \to \mathrm{U}^{(k)}$ in probability 
as $T\to \infty$, so it must be that $\mathrm{U}^{(k)}\eqdist \tilde{\mathrm{U}}^{(k)}$. Thus, we conclude that $\mathrm{U}^{(n,k)}\convdist \mathrm{U}^{(k)}$ in the product topology as $n\to \infty$, proving \eqref{conv of weights w}.
\end{proof}

\begin{lemma}\label{rhino appl}
For $k\ge 1$, it holds that 
\begin{equation*}
\lim_{T\to \infty}\liminf_{n\to\infty}\Pp(\hvcls{j}_T  = \hvcls{j}\text{ for all $j\in [k]$}) =1.
\end{equation*}
\end{lemma}
\begin{proof}
We use the following deterministic fact. Suppose that $(x^{(n,i)}_T,i \in [k])_{n \ge 1,T>0}$, $(x^{(n,i)},i \in [k])_{n \ge 1})$, $(x^{(i)}_T,i \in [k])_{T > 0}$, and $(x^{(i)},i \in [k])$ are decreasing vectors of real numbers such that 
\begin{enumerate}
\item $(x^{(i)},i \in [k])$ is strictly decreasing;
\item $(x^{(i)}_T,i \in [k])=(x^{(i)},i \in [k])$ for all $T$ sufficiently large;
\item for all $T>0$, $(x^{(n,i)}_T,i \in [k])_{n \ge 1}$ converges to $(x^{(i)}_T,i \in [k])$ pointwise as $n \to \infty$;
\item $(x^{(n,i)},i \in [k])_{n \ge 1}$ converges to $(x^{(i)},i \in [k])$ pointwise as $n \to \infty$; and
\item for all $n \ge 1$, $T>0$, and $i \in [k-1]$, if $x^{(n,i)}_T \ne x^{(n,i)}$ then $x^{(n,i)}_T \le x^{(n,i+1)}$.
\end{enumerate}
Then $(x^{(n,i)}_T,i \in [k-1])=(x^{(n,i)},i \in [k-1])$ for all $T$ and $n$ sufficiently large. We explain why this is true in the case $k=2$, as the logic for general $k$ is very similar. Write $\veps=x^{(1)}-x^{(2)}>0$, and fix $T_0$ large enough that $(x^{(1)}_T,x^{(2)}_T)=(x^{(1)},x^{(2)})$ for all $T>T_0$. Then for $T>T_0$ we have that for all $n$ sufficiently large, $|x^{(n,i)}_T-x^{(i)}|<\tfrac{\veps}{4}$ and $|x^{(n,i)}-x^{(i)}|<\tfrac{\veps}{4}$ for $i \in \{1,2\}$. Now, if $x^{(n,1)}_T \ne x^{(n,1)}$ then 
\[
x^{(n,1)}_T \le x^{(n,2)} < x^{(2)}+\frac{\veps}{4} \le x^{(1)}-\frac{3\veps}{4} < x^{(n,1)}_T-\frac{\veps}{2}\, ,
\]
a contradiction; so $x^{(n,1)}_T=x^{(n,1)}$. 

We apply the deterministic fact with $x^{(i)}=\excsize{i}$, $x^{(i)}_T=\excsize{i}_T$, $x^{(n,i)}=w(\hvcls{i})$, and $x^{(n,i)}_T=w(\hvcls{i}_T)$. The first condition in the above list holds since $(\excsize{i},i \ge 1)$ is almost surely strictly decreasing. The second condition holds by Lemma~\ref{brownian excs}. By the Skorohod representation theorem and Proposition~\ref{joint convergence to exc kT}, we may work in a space in which the third condition holds almost surely, and, likewise, by~Lemma~\ref{new aldous} we may work in a space in which the fourth condition also holds almost surely. The fifth condition is automatic from the fact that $(\hvcls{i},i\ge 1)$ lists the components of $\hat{G}_n$ by weight, and $(\hvcls{i}_T,i\ge 1)$ lists a subset of the same set of components by weight. 

This shows that 
\[
\lim_{T\to \infty}\liminf_{n\to\infty}\Pp(w(\hvcls{j}_T)  = w(\hvcls{j})\text{ for all $j\in [k]$}) =1.
\]
The lemma follows from this and the fact that $w(\hvcls{j+1})<w(\hvcls{j})$ for all $j \in [k-1]$ with probability tending to $1$ as $n \to \infty$.
\end{proof}

Recalling the definitions of $\clvec_T$ and $\clvec$ from~\eqref{def: W-vector}, \cref{rhino} and \cref{rhino appl} together imply the following corollary.
\begin{cor}\label{remove T whp}
For all $k\ge 1$, it holds that 
\[
\lim_{T\to \infty}\liminf_{n\to\infty}\Pp(\clvec_T = \clvec) = 1.
\]
\end{cor}
We are now ready to prove the convergence~\eqref{eqcv: main} in the main result, albeit in a weaker topology. 
\begin{lemma}\label{pt convergence}
For all $k\ge 1$, it holds that $\clvec \convdist \exvec$ as $n\to \infty$. 
\end{lemma}
\begin{proof}
Fix $k\geq 1$. 
Recall from \cref{conv of components} that $\clvec_T \convdist \exvec_T$ as $n\to \infty$. 
It follows from \cref{remove T whp} that
\[
\lim_{T\to \infty}\limsup_{n\to\infty}
\Pp\left( \clvec_T\ne  \clvec \right)=0\, ,
\]
so by the principle of accompanying laws, there is a random variable $\limitvec$ such that $\exvec_T \to \limitvec$ weakly as $T\to\infty$, and $\clvec \to \limitvec$ in distribution as $n\to\infty$. Since $\exvec_T \to \exvec$ in probability as $T\to \infty$ by \cref{brownian excs}, we conclude that  $\clvec \to \exvec$ in distribution as $n\to\infty$.
\end{proof}

\subsection{$\ell^2$-convergence}
\label{sec: l2}
Having already obtained the convergence in the product topology, in this section we prove the full statement of~\cref{samelaw}.
Our proof builds on the idea of Aldous~\cite{Aldous97}, where we construct a size-biased point process (SBPP) on $\R_+ \times (0,\infty)$ from the sequence of component weights and argue that this SBPP converges. We then use \cite[Proposition 15]{Aldous97}, which is restated below as Theorem~\ref{thm: sbpp l2}, to conclude the proof. Before we prove the $\ell^2$-convergence result, we introduce the notion of convergence we require for the SBPP, and the construction of the SBPP.

Let $\mathcal{M}$ be the space of Radon measures on $\R{_+} \times (0,\infty)$, equipped with the topology of vague convergence. 
For a random point process $X$ on $\R\times(0,\infty)$, write $\mu(X)=\sum_{(s,x) \in X} \delta_{(s,x)}$ for the associated ($\mathcal{M}$-valued) random counting measure. Then for a collection of random point processes $(X_n,1 \le n \le \infty)$ on $\R_+ \times (0,\infty)$, we say that $X_n \to X_\infty$ in distribution if  $\mu(X_n) \to \mu(X_{\infty})$ in distribution. 

Given a random $\ell^2$-valued sequence $X = (X_i)_{i\in \N}$ of non-negative real numbers, we define an associated SBPP as follows: conditional on $(X_i)_{i\in \N}$, let $(E_i)_{i\in \N}$ be a collection of independent exponential variables of respective rates $\left(X_i, i\in \N\right)$. Then, for each $i\in \N$, set 
\[
S_i =\sum_{j\in \N}X_j\mathbf 1_{\{E_j\le E_i\}},
\]
which is finite almost surely since $X \in \ell^2$. Define 
\[
\Pi_X = \{(S_i, X_i): i\in \N\}.
\]
We shall call this resulting point process $\Pi_X$ the {\em size-biased point process generated by the sequence $X$}. Then, $\Pi_X$ is a random point measure on $\R_+ \times(0,\infty)$ that a.s.\ satisfies the following properties (see \cite[Section 3.3]{Aldous97}):

\begin{enumerate}[(i)]
\item 
for any $t\ge 0$ and $\delta\in(0,1)$, $\big|\Pi_X\cap [0, t]\times [\delta,1/\delta]\big| <\infty$ almost surely;
\item
for each $(s, x)\in \Pi_X$, we have
\begin{equation}
\label{cond1: SBPP}
\sum_{(s',x')\in\Pi_X}x'\mathbf 1_{\{s'\le s\}}= s \quad \text{a.s.};
\end{equation}
\item 
we have
\begin{equation}
\label{cond2: SBPP}
\sup\{x: (s,x)\in \Pi_X \text{ for some }s>s_0\} \xrightarrow{s_0\to\infty} 0 \quad \text{in probability}.
\end{equation}
\end{enumerate}

Note that, given an SBPP $\Pi_X$ associated to a random sequence $X=(X_i)_{i\in \N}$ of non-negative real numbers in $\ell^2$, we can obtain a sequence of the $X_i$'s in a decreasing order, by ranking the second coordinates of atoms of $\Pi_X$. Formally, we denote by $\mathrm{ord}(\Pi_X)$ this ranked sequence.  Then, $\mathrm{ord}(\Pi_X)$ is an element of $\ell^2_{\searrow}\coloneqq\{(x_i)_{i\in \N}: x_1\ge x_2\ge x_3\ge \cdots\ge 0, \sum_{i\in \N}x_i^2<\infty\}$, the space of $\ell^2$-sequences arranged in decreasing order.  

We will use the following result from Aldous~\cite{Aldous97} to show the $\ell^2$-convergence of the component sizes using vague convergence of a related SBPP.

\begin{thm}[Proposition 15,~\cite{Aldous97}] \label{thm: sbpp l2}
For each $n\in \N$, let $X^{(n)}=(X_{n, i})_{i\in\N}$ be a random sequence taking values in $\ell^2$ and let $\Pi_n := \Pi_{X^{(n)}}$ be the associated size-biased point process. Let $\Pi$ be a random point process on $\R_+ \times (0,\infty)$ satisfying~\eqref{cond1: SBPP},~\eqref{cond2: SBPP} and that
\begin{equation} \label{cond3: SBPP}
\sup\{s \in \R_+:(s,x)\in \Pi \text{ for some } x\}=\infty \quad \text{a.s.}
\end{equation}
Suppose that $\Pi_n\to\Pi$ in distribution. Then $\mathrm{ord}(\Pi)$ takes values in $\ell^2_{\searrow}$ almost surely and $\mathrm{ord}(\Pi_n)\to \mathrm{ord}(\Pi)$ in distribution with respect to the $\ell^2$-topology. 
\end{thm}
Recall from~\eqref{explore-times} that component sizes in the size-biased exploration can be expressed in terms of $\sscount{t}$ and $\hat{\wst}_k$. 
With the understanding that $\sscount{0-}=0$, we define
\[
\Pi_n=\Big\{\Big(n^{-2/3}\sscount{\hat{\wst}_{k}-}, n^{-2/3}\big(\sscount{\hat{\wst}_{k+1}-}-\sscount{\hat{\wst}_k-}\big)\Big): 0\le k< \hat{\wK}_n-1\Big\}\, .
\]
As the components appear in a size-biased manner, $\Pi_n$ is a size-biased point process in the above sense. 
In particular, ranking the second coordinates of $\Pi_n$ in decreasing order yields the sequence $(n^{-2/3}|\lgcls{i}|: i\ge 1)$. 
Recalling the notation $g^{(i)}, \excsize{i}$, we denote 
\[
\Pi_{\infty}=\big\{\big(g^{(i)}, \excsize{i}\big): i\ge 1\big\},
\]
which almost surely has only finitely many atoms on $[0, t]\times [\veps, \infty)$ for any $t, \veps>0$. Further, $\Pi_{\infty}$ satisfies~\eqref{cond1: SBPP},~\eqref{cond2: SBPP} and~\eqref{cond3: SBPP} (see Section 3.4 in~\cite{Aldous97}).

We claim the following 
\begin{lemma}
\label{lem: cv-Pi}
It holds that
\[
\Pi_n \convdistn \Pi_{\infty}.
\]
\end{lemma}     
\begin{proof}[Proof of Lemma~\ref{lem: cv-Pi}]
Since almost surely $\Pi_{\infty}$ is a simple point measure with no fixed point of jump, the claimed convergence will follow once we show that for all $T,\veps>0$, 
\begin{equation}
\label{eq: cv-Pin}
\Pi_n\cap \big([0, T]\times [\veps, \infty)\big) \convdistn \Pi_{\infty}\cap \big([0, T]\times [\veps, \infty)\big)\, ;
\end{equation}
see \cite[Propositions~3.12 and~3.13]{Resnick1987} together with \cite[Theorem~2.7]{Loehr2013}.

For $k\in \N$ and $\veps>0$, denote by $\mathcal A(k, \veps)$ the event that $\excsize{k+1}\ge \veps$. As the sequence $(\excsize{i}, i\ge 1)$ converges to 0 almost surely, we can choose $k_0=k_0(\veps)$ such that $\Pp(\mathcal A(k_0, \veps))\le \veps$. Recall that $\exc{i}_T=(g^{(i)}_T, d^{(i)}_T)$ stands for the $i$th longest excursion of $\bm$ that has ended by time $T$. 
On the event $\mathcal{A}(k_0, \veps)^c$, we can write 
\[
\Pi_{\infty}\cap\big([0, T]\times [\veps, \infty)\big) = \big\{\big(g^{(i)}_T, \excsize{i}_T\big): i\in [k'_0]\big\},
\]
with $k'_0\le k_0$. 

Next, recall that $\hvcls{i}_T$ is the $i$th heaviest connected component explored by time $Tn^{2/3}$ in the size-biased exploration, and that its exploration is carried out during the interval $(\hat g^{(n, i)}_T, \hat d^{(n, i)}_T)$. Also recall that $\lgcls{i}_T$ is the $i$th largest connected component explored by time $Tn^{2/3}$ in the size-biased exploration. Let $\mathcal A_n(k_0)$ be the event that there is an  $i\in [k_0]$ such that $\hvcls{i}_T\ne \lgcls{i}_T$. By Lemma \ref{hvT equal lgT} we may choose $n$ large enough so that $\Pp(\mathcal A_n(k_0))<\veps$.

Thanks to~\cref{joint convergence to exc kT}, we have
\[
\Big\{n^{-2/3}\Big(\hat g^{(n, i)}_T, (\hat d^{\,(n, i)}_T-\hat g^{(n,i)}_T)\Big), i\in [k_0]\Big\} \convdistn \Big\{\big(g^{(i)}_T, \excsize{i}_T\big), i\in [k_0]\Big\}.
\]
Combining this with~\cref{LLN 2}, we further deduce that 
\[
\Big\{n^{-2/3}\Big(\hat N(\hat g^{(n, i)}_T-), \big(\hat N(\hat d^{\,(n, i)}_T-)-\hat N(\hat g^{(n,i)}_T-)\big)\Big), i\in [k_0]\Big\} \convdistn \Big\{\big(g^{(i)}_T, \excsize{i}_T\big), i\in [k_0]\Big\}. 
\]
In particular, the convergence still holds when considering the restriction of both sides on $[0, T]\times [\veps, \infty)$. This proves the convergence in~\eqref{eq: cv-Pin} on $\left(\mathcal A(k_0, \veps)\cup \mathcal A_n(k_0)\right)^c$. The proof is now complete by noting that $\Pp(\mathcal A(k_0, \veps)\cup \mathcal{A}_n(k_0))\le 2\veps$ for $n$ large enough.
\end{proof}

\begin{remark}
\label{rem: SBPP}
Let us compare our construction of $\Pi_n$ with the ones in Aldous~\cite{Aldous97} and in Bhamidi,Van der Hofstad and Van Leeuwaarden~\cite{BHL2010}. In~\cite{Aldous97}, the following point process was considered instead:
\[
\hat\Pi_n = \{n^{-2/3}\tau_k, n^{-2/3}(\tau_{k+1}-\tau_k): 0\le k<K_n-1\}.
\]
Since $\tau_{k+1}-\tau_k$ corresponds to the weight of the $k$th component explored in the weight-biased exploration, $\hat\Pi_n$ is also a size-biased point process as defined previously. By establishing the convergence of $\hat\Pi_n$ to $\Pi_{\infty}$, Aldous thus proved the $\ell^2$-convergence of the component weights. In~\cite{BHL2010}, the authors introduced a discrete-time exploration process where excursion lengths coincided with the component sizes. The main issue in the application of the size-biased point process technique in~\cite{BHL2010} stems from the fact that in the exploration considered there, the components appear in a weight-biased order, rather than in a size-biased order. In \cite{BHL2010}, the proofs for convergence of the rescaled component sizes in both the $\ell^2$-topology \emph{and} the product topology rely on this erroneous application of the point process technique. 
\end{remark}

\begin{proof}[Proof of \cref{samelaw}]
As an immediate consequence of Lemma~\ref{lem: cv-Pi} and Theorem~\ref{thm: sbpp l2} and the fact that $(|\lgcls{i}|, i\ge 1)\eqdist(|\lgclw{i}|, i\ge 1)$, we have
\[
\left(n^{-2/3}|\lgclw{i}|, i\ge 1\right) \convdistn \big(\excsize{i},i\ge 1\big), 
\]
with respect to the $\ell^2$-topology. From \cite[Proposition 4]{Aldous97}, 
we then get the convergence of $(n^{-2/3}w(\hvclw{i}),i\ge 1)$ to $(\excsize{i},i\ge 1)$ in distribution in the $\ell^2$-topology. 

We now claim that $(n^{-2/3}|\hvclw{i}|,i\ge 1)$ and $(n^{-2/3}w(\lgclw{i}), i\ge 1)$ also jointly converge to $(\excsize{i},i\ge 1)$ in distribution in the $\ell^2$-topology. 

We first show how the $\ell^2$-convergence of $(n^{-2/3}|\lgclw{i}|,i\ge 1)$ implies the $\ell^2$-convergence of $(n^{-2/3}|\hvclw{i}|,i\ge 1)$. Fix $\veps>0$. 
Since $(n^{-2/3}|\lgclw{i}|,i\ge 1)$ converges in $\ell^2$, in particular the sequence is tight. As a result, we can pick $k_0$ large enough that
\begin{equation}
\label{eq_ell2_10}
\Pp\Big(\sum_{i>k_0} |\lgclw{i}|^2>\veps n^{4/3}\Big)<\veps.
\end{equation}
Let $A_{n, k_0}$ be the event that $\lgclw{i}=\hvclw{i}$ for each $i\in [k_0]$. By \cref{cor_consist}, which we recall we are allowed to use since we already proved that the convergence \eqref{eqcv: main} holds in the product topology, we have that $\Pp(A_{n, k_0})\ge 1-\veps$ for $n$ large enough. 
Since $(|\hvclw{i}|,i\ge 1)$ is a rearrangement of $(|\lgclw{i}|,i\ge 1)$, on the event $A_{n, k_0}$, we have
\[
\sum_{i>k_0}|\hvclw{i}|^2 = \sum_{i\ge 1}|\hvclw{i}|^2-\sum_{i\in [k_0]}|\hvclw{i}|^2=\sum_{i\ge 1}|\lgclw{i}|^2-\sum_{i\in [k_0]}|\lgclw{i}|^2=\sum_{i>k_0}|\lgclw{i}|^2. 
\]
It then follows from~\eqref{eq_ell2_10} that for $n$ large enough,
\[
\Pp\Big(\sum_{i>k_0} |\hvclw{i}|^2>\veps n^{4/3}\Big)<2\veps,
\]
which shows that $(n^{-2/3}|\hvclw{i}|,i\ge 1)$ is also tight in $\ell^2$. Lemma~\ref{pt convergence} implies that the only limit point is $(\excsize{i}, i\ge 1)$. We conclude that $(n^{-2/3}|\hvclw{i}|,i\ge 1)$ converges in distribution to $(\excsize{i},i\ge 1)$ with respect to the $\ell^2$-topology. 

The argument that the $\ell^2$-convergence of $(n^{-2/3}w(\hvclw{i}), i\ge 1
)$ implies the $\ell^2$-convergence of $(n^{-2/3}w(\lgclw{i}), i\ge 1)$ is analogous, so we omit it here.

We have shown that each component in the left-hand side of~\eqref{eqcv: main} individually converges in distribution to $(\excsize{i}, i\ge 1)$ with respect to the $\ell^2$-topology. In particular, this implies that the left-hand side of~\eqref{eqcv: main} is tight in $\ell^2$.
However, \cref{pt convergence} implies that the limit must be the same for all four coordinates, so the joint convergence in~\eqref{eqcv: main} indeed holds.
\end{proof}

\section{Studying the root weights}\label{root wt proofs}
In this section, we study the processes that encode the weights of the roots when the roots are selected either uniformly at random or in a weight-biased manner to prove~\cref{roots eps dense}, \cref{roots dont matter} and~\cref{prop: comps explored}. For the sake of clarity and to prevent circular reasoning, we emphasise that the only result from Section~\ref{Sub:product topology} that we rely on here is Proposition~\ref{prop:conv_S_n}.

Recall that $\wroot_1,\wroot_2,\dots$ is the sequence of roots when they are selected in a weight-biased manner from the unexplored vertices and $\sroot_1,\sroot_2,\dots$ is the sequence of roots when they are selected uniformly at random from the unexplored vertices. The next proposition states that in both cases, the cumulative root weight process converges under rescaling. 

\begin{prop}\label{prop:roots_uar}
    It holds that
    \[\left(\frac{1}{n^{1/3}}\!\sum_{i=1}^{\lfloor\!sn^{1/3}\!\rfloor}\!\wt{\wroot_i}, s\ge 0\right) \convprobn  (s, s\ge 0)\text{ and }\!\left(\frac{1}{n^{1/3}}\!\sum_{i=1}^{\lfloor\!sn^{1/3}\!\rfloor}\!\wt{\sroot_i}, s\ge 0\right) \convprobn  (\mu s, s\ge 0)\]
    uniformly on compact sets. 
\end{prop}

We will use the next proposition to show that, for both root selection procedures, both the number of root vertices that contribute to $S_n$ and the sum of their weights are negligible on our scale of interest.

\begin{prop}\label{prop:corr_uar}
    For any $s\ge 0$ and $t \ge 0$, 
    \begin{align*}
    &n^{-1/3}\sum_{i=1}^{\lfloor sn^{1/3}\rfloor} \wt{\wroot_i}\I{J_{\wroot_i }\le t n^{2/3}}, \quad n^{-1/3}\sum_{i=1}^{\lfloor sn^{1/3}\rfloor} \wt{\sroot_i}\I{J_{\sroot_i} \le t n^{2/3}} \\
    & n^{-1/3}\sum_{i=1}^{\lfloor sn^{1/3}\rfloor} \I{J_{\wroot_i }\le t n^{2/3}}, \quad \text{and}\quad n^{-1/3}\sum_{i=1}^{\lfloor sn^{1/3}\rfloor} \I{J_{\sroot_i} \le t n^{2/3}}
    \end{align*}
    all converge to 0 in probability as $n\to\infty$.  
\end{prop}

Let us point out that, together with Proposition~\ref{prop: comps explored}, the previous two propositions immediately yield  Lemmas~\ref{roots eps dense} and~\ref{roots dont matter}.

\begin{proof}[Proof of Lemma~\ref{roots eps dense}]
This follows from~\cref{prop:roots_uar} and~\cref{prop: comps explored}.
\end{proof}
\begin{proof}[Proof of Lemma~\ref{roots dont matter}]
This follows from~\cref{prop:corr_uar} and~\cref{prop: comps explored}.
\end{proof}

To prove Propositions~\ref{prop:roots_uar} and~\ref{prop:corr_uar}, we will sample the order of the root vertices as follows. For each $i\in [n]$ we sample  $E_i\sim \operatorname{Exp}(\wt{i}/\sum_{j=1}^n \wt{j})$ and $\hat{E}_i\sim \operatorname{Exp}(1/n)$.
By the properties of competing exponentials, if we consider the permutation $\Sigma$ of $[n]$ that satisfies that $(E_{\Sigma(k)},k=1,\dots n)$ is increasing, then $\wtv_{\Sigma}=(\wt{\Sigma(1)},\dots, \wt{\Sigma(n)})$ is a size-biased random reordering of $\wtv=(\wt{1},\dots,\wt{n})$. Similarly, if we consider the permutation $\hat{\Sigma}$ of $[n]$ that satisfies that $(\hat{E}_{\hat{\Sigma}(k)},k=1,\dots n)$ is increasing then $(\wt{\hat{\Sigma}(1)},\dots, \wt{\hat{\Sigma}(n)})$ is a uniformly random reordering of $(\wt{1},\dots,\wt{n})$. 

We may then define $(R_1,\ldots,R_{K_n})$ as follows, Set $\wroot_1=\Sigma(1)$.  Then, recursively, conditional on $\wroot_1,\dots, \wroot_j$ and $J_1,\dots, J_n$ we let $\wroot_{j+1}$ be the first vertex that occurs after $\wroot_j$ in $\Sigma$ that has not been included in the exploration by time $\wst_j$. Formally,
\[\wroot_{j+1}=\Sigma(\min\{k>\Sigma^{-1}(\wroot_{j}):J_k>\wst_j \} ).\]
Similarly, we may define $(\hat{R}_1,\ldots,\hat{R}_{\hat{K}_n})$ by setting
$\sroot_1=\hat{\Sigma}(1)$ and, for $j\ge 1$, conditional on $\sroot_1,\dots, \sroot_j$ and $J_1,\dots, J_n$, setting 
\[\sroot_{j+1}=\hat{\Sigma}(\min\{k>\hat{\Sigma}^{-1}(\sroot_{j}):J_k>\sst_j \}).\]

Let $V(t) \coloneqq \sum_{i=1}^n\wt{i}\I{E_i\le t }$ record the total weights of the vertices $i$ whose exponential clock $E_i$  rings before time $t$. Let $T_k \coloneqq \inf\{t \ge 0 :\left|\{i:E_i\le t\}\right|\ge k\}$ be the time that the $k$th exponential clock rings; in other words, $T_k=E_{\Sigma(k)}$.  Similarly, define $\hat{V}(t) \coloneqq \sum_{i=1}^n\wt{i}\I{\hat{E}_i\le t }$ and $\hat{T}_k \coloneqq \inf\{t \ge 0:|\{i:\hat{E}_i\le t\}|\ge k\}=E_{\hat{\Sigma}(k)}$. 

\begin{lemma}\label{lem:lln_roots}
It holds that 
\begin{align}
\label{Vn conv}    \left( n^{-1/3}V(n^{1/3}s), s\ge 0\right)& \convprobn (s, s\ge 0),\\
    \label{Vn conv 2}\left( n^{-1/3}\hat{V}(n^{1/3}s), s\ge 0\right)& \convprobn ( \mi s, s\ge 0),\\ \label{Tn conv}
    \left( n^{-1/3}T_{\lfloor n^{1/3}s\rfloor }, s\ge 0\right)& \convprobn ( s, s\ge 0),\text{ and}\\ \label{Tn conv2}
     \left( n^{-1/3}\hat{T}_{\lfloor n^{1/3}s\rfloor }, s\ge 0\right)& \convprobn ( s, s\ge 0)
\end{align}
uniformly on compact sets.
\end{lemma}

The proof of Lemma~\ref{lem:lln_roots} appears in Appendix~\ref{sec: comment}. 
With this lemma in hand, to deduce Proposition \ref{prop:roots_uar} we need to show that correcting jumps in $V$ and $\hat{V}$ due to vertices that have already caused a jump in $S_n$ has a negligible effect on the time horizon and scale of interest. Similarly, to deduce Proposition \ref{prop:corr_uar} we need to show that correcting jumps in $S_n$ coming from vertices that have already been used as a root vertex has a negligible effect on the time horizon and scale of interest. For both, we need the next lemma. 
\begin{lemma}\label{lem:competing_clocks}
For any $s\ge 0$ and $t\ge 0$, the following sums 
\begin{align*}
&n^{-1/3}\sum_{i\in[n]} \wt{i} \I{J_i\le tn^{2/3}, E_i\le sn^{1/3} }, \quad 
n^{-1/3}\sum_{i\in[n]} \wt{i} \I{J_i\le tn^{2/3}, \hat{E}_i\le sn^{1/3} },\\
&n^{-1/3}\sum_{i\in[n]} \I{J_i\le tn^{2/3}, E_i\le sn^{1/3} } ,\text{ and }
n^{-1/3}\sum_{i\in[n]}  \I{J_i\le tn^{2/3}, \hat{E}_i\le sn^{1/3} }
\end{align*}
all converge to 0 in probability as $n\to\infty$. 
\end{lemma}
\begin{proof}
This follows immediately from the first moment method. Fix $s$ and $t$.  Then, using that $1-\exp(-x)\le x$ for all $x\ge 0$, we see that
\begin{align*}
&\quad\mathbb E\left[\sum_{i\in[n]} \wt{i} \I{J_i\le tn^{2/3}, \hat{E}_i\le sn^{1/3} }\right]\\
&=\sum_{i\in[n]} \wt{i} \left(1-\exp\left(-t n^{2/3} \wt{i} \frac{(1+\cw n^{-1/3})}{\sum_{j\in[n]} \wt{j}}\right)\right)(1-\exp(-sn^{-2/3}))\\
&\le \frac{ts(1+|\cw| n^{-1/3})\sum_{i\in[n]} (\wt{i})^2}{\sum_{i\in[n]} \wt{i}}=O(1),
\end{align*}
by the assumptions \eqref{first moment} and \eqref{second moment}. Therefore, 
\[\mathbb E\left[n^{-1/3}\sum_{i\in[n]} \wt{i} \I{J_i\le tn^{2/3}, \hat{E}_i\le sn^{1/3} } \right]\to 0\]
as $n\to \infty$ and the first claim follows.

Similarly, 
\[\mathbb E\left[\sum_{i\in[n]} \wt{i} \I{J_i\le tn^{2/3}, E_i\le sn^{1/3} }\right]\le \frac{nts(1+|\cw| n^{-1/3})\sum_{i\in[n]} (\wt{i})^3}{(\sum_{i\in[n]} \wt{i})^2}=O(1)
\]
by assumptions \eqref{first moment} and \eqref{third moment}, which implies the second claim.
The third and fourth claims follow analogously. 
\end{proof}

We now combine Lemmas \ref{lem:lln_roots} and \ref{lem:competing_clocks} to prove Proposition \ref{prop:roots_uar}. 
\begin{proof}[Proof of Proposition~\ref{prop:roots_uar}]
    Fix $t>0$. Recall that $\hat{V}(\hat{T}_k)$ is the sum of the first $k$ increments of $\hat{V}$. To prove the second assertion of the proposition, we will show that
    \begin{equation}\label{eq:rootsumvdiff}
        \sup_{0\le k \le\lfloor tn^{1/3}\rfloor } n^{-1/3}\left|\sum_{i\in[k]} \wt{\sroot_i} - \hat{V}(\hat{T}_k)\right|\convprobn 0\, ;
    \end{equation}
    the second assertion then follows due to the convergence under rescaling of $\hat{V}$ from Lemma \ref{lem:lln_roots}.
    
    We recall the difference between the processes $\sum_{i\in[k]} \wt{\sroot_i}$ and $\hat{V}(\hat{T}_k)$. In the process $\sum_{i\in[k]} \wt{\sroot_i}$, it is taken into account that the $j$th root \emph{cannot} be a vertex that has already appeared in the exploration by time $\wst_j$, whereas the process $\hat{V}(\hat{T}_k)$ does not reject such vertices. We do, however, obtain the next increment of $\sum_{i\in[k]} \wt{\sroot_i}$ by taking the next increment of $\hat{V}$ that we do not reject. Thus, to show that $\sum_{i\in[k]} \wt{\sroot_i}$ and $\hat{V}(\hat{T}_k)$ are close, we need to show that not too many increments get rejected and that their sum is not too big. 

    More concretely, the main fact that we use is as follows. Fix $k, r\in \N$ and $x, y, z\in (0,\infty)$, and let $E=E(k,r,x,y,z)$ be the event that:
    \begin{enumerate}
    \item\label{enum_rt_1} $\inf_{0\le t\le x}S_n(t)\le -y $
     \item \label{enum_rt_4} $\hat{V}(\hat{T}_{k+r})\le y$
     \item \label{enum_rt_3} $\hat{T}_{k+r}\le z$
    \item\label{enum_rt_2} $\sum_{i\in[n]} \I{J_i\le x , \hat{E}_i\le z}\le r$.
    \end{enumerate}
    Then, on the event $E$, 
     it holds that
    $\hat{E}_{\hat{R}_k}\le z$ 
    and that for all $\ell \in [k]$,
    \begin{equation}\label{bound_rt_w}\hat{V}(\hat{T}_\ell)- \sum_{i\in[n]} \wt{i}\I{J_i\le x , \hat{E}_i\le z} \le \sum_{i\in[\ell]} \wt{\sroot_i}\le \hat{V}(\hat{T}_{\ell+r}).\end{equation}
    Let us explain why this is true.
    Statements (\ref{enum_rt_1}) and (\ref{enum_rt_4}) together imply that, if no more than $r$ roots have been rejected for the first $k$ components, then the $k$th component is fully explored by time $x$. Indeed, if at most $r$ roots have been rejected, then $\hat{V}(\hat{T}_{k+r})$ is an upper bound for the cumulative root weights of the first $k$ components, and $\hat{V}(\hat{T}_{k+r})$ is in turn bounded from above by $y$.     
    Recall that in the size-biased exploration process, the exploration of the first $k$ components ends at time $\hat{\wst}_k$, which is the first moment when $S_n$ drops below the level $-\sum_{j\in [k]}\wt{\sroot_j}+\sum_{j\in [k]} \wt{\sroot_j}\I{J_{\sroot_j}\le \hat{\wst}_k}$.
    Thus, when (1) and (2) occur, the exploration of the first $k$ components ends before $S_n$ reaches level $-y$, which happens before time $x$.
    Next, (\ref{enum_rt_3}) implies that the first $k+r$ jumps of $\hat{V}$ occur before time $z$. Finally, (\ref{enum_rt_2}) implies that no more than $r$ roots get rejected before time $x$ in $S_n$ and time $z$ in $\hat{V}$. Thus, we conclude that no more than $r$ roots get rejected for the first $k$ components and that the $k$th component is explored before time $x$ in $S_n$ and time $z$ in $\hat{V}$,  and in particular that $\hat{E}_{\hat{R}_k}\le z$. Thus, for $\ell \le k$, the sum $\hat{V}(\hat{T}_{\ell+r})$ of the first $\ell+r$ jumps of $\hat{V}$ is an upper bound for the sum of the first $\ell$ root weights, which yields the second bound of \eqref{bound_rt_w}. Moreover, the sum of the first $\ell$ jumps $\hat{V}(\hat{T}_\ell)$, minus the weights of all rejected jumps by time $x$ in $S_n$ and time $z$ in $\hat{V}$, is a lower bound for the sum of the first $\ell$ root weights, which yields the first bound of \eqref{bound_rt_w}.
    
    Fix $\veps>0$. We will now choose $r,x,y,z$ such that, for $n$ large enough, (\ref{enum_rt_1})--(\ref{enum_rt_2}) occur with probability at least $1-\veps$ and such that the upper and lower bound in \eqref{bound_rt_w} differ by at most $\veps n^{1/3}$ for all $1\le \ell \le \lfloor tn^{1/3}\rfloor$ with probability at least $1-\veps$. 
    
    First, by the convergence under rescaling of $\hat{V}$ and $\hat{T}$ in Lemma \ref{lem:lln_roots},  we may pick $0<\delta<t$ small enough and $n$ large enough that 
    \begin{equation}\label{eq:rw_1}\p{\sup_{\ser{1}\le \ell \le \lfloor tn^{1/3}\rfloor} \left(\hat{V}(\hat{T}_{\ell+\lfloor \delta n^{1/3}\rfloor })-\hat{V}(\hat{T}_{\ell})\right)>\veps n^{1/3}}<\frac{\veps}{2}.
    \end{equation}

    Since $W(u)\to -\infty$ almost surely as $u\to \infty$, we may pick $t$  large enough that \[\p{\inf_{0\le u\le t} W(u)>-4\mi t }<\frac{\veps}{8}.\]

    By the uniform convergence in distribution of $\scaledrw$ to $W$ in Proposition~\ref{prop:conv_S_n} and the choice for $t$, we may pick $n$ large enough such that
    \begin{equation}\label{eq:rw_2}\p{\inf_{0\le u\le t} \scaledrw(u)>-3\mi s }<\frac{\veps}{4} .\end{equation}

    By Lemma~\ref{lem:competing_clocks} we may pick $n$ large enough such that 
    \begin{equation}\label{eq:rw_3}\p{\sum_{i\in[n]} \wt{i}\I{J_i\le tn^{2/3} , \hat{E}_i\le 2sn^{1/3}}>\veps n^{1/3}}<\frac{\veps}{2}.\end{equation}
    
    and
    \begin{equation}\label{eq:rw_4}\p{\sum_{i\in[n]}\I{J_i\le tn^{2/3},\hat{E}_i\le 2sn^{1/3}}> \lfloor \delta n^{1/3}\rfloor }<\frac{\veps}{4}.\end{equation}
   
    Moreover, since we chose $\delta<t$, by Lemma \ref{lem:lln_roots} we may pick $n$ large enough such that 
    \begin{equation}\label{eq:rw_5}\p{\hat{T}_{\lfloor sn^{1/3} \rfloor +\lfloor \delta n^{1/3}\rfloor}>2sn^{1/3}}\le \frac{\veps}{4}.\end{equation}
    
    By Lemma \ref{lem:lln_roots}, we also know that 
    \[ 
    n^{-1/3} \hat{V}(\hat{T}_{\lfloor (s+\delta)n^{1/3}\rfloor}) \convprobn \mi \cdot (s+\delta),
    \]
    so we can make $n$ large enough that 
    \begin{equation}\label{eq:rw_6}\p{\hat{V}(\hat{T}_{\lfloor (s+\delta) n^{1/3}\rfloor})> 3\mi sn^{1/3}}<\frac{\veps}{4}.\end{equation}
  Choosing $r=\lfloor \delta n^{1/3}\rfloor$, $k=\lfloor sn^{1/3}\rfloor$, $x=tn^{2/3}$, $y=3n^{1/3}\mi s$ and $z=2sn^{1/3}$, it follows from~\eqref{eq:rw_1}--\eqref{eq:rw_6} that (\ref{enum_rt_1})--(\ref{enum_rt_2}) hold with probability at least $1-\veps$ and that the upper and lower bound in \eqref{bound_rt_w} differ by at most $\veps n^{1/3}$ for all $0\le \ell \le \lfloor sn^{1/3}\rfloor$ with probability at least $1-\veps$. 
  Since $\veps>0$ and $\delta>0$ were arbitrary, this establishes \eqref{eq:rootsumvdiff} and thus the second statement in Proposition \ref{prop:roots_uar}. The proof for the first statement is completely analogous.     
\end{proof}
Lemma~\ref{lem:lln_roots} and Lemma~\ref{lem:competing_clocks} now straightforwardly yield Proposition \ref{prop:corr_uar}.

\begin{proof}[Proof of Proposition~\ref{prop:corr_uar}]
To prove the second convergence, let us denote $\sroot_{\ast}=\sroot_{\lfloor sn^{1/3}\rfloor}$ and observe that for any $\veps>0$,
\begin{align*} 
\Pp\left(n^{-1/3}\sum_{i=1}^{\lfloor s n^{1/3}  \rfloor} \wt{\sroot_i}\I{J_{\sroot_i} \le tn^{2/3}}>\veps \right)
&\le \Pp\left(n^{-1/3}\sum_{i=1}^{n} \wt{i}\I{J_i \le tn^{2/3}, \hat{E}_i \le \hat{E}_{\sroot_{\ast}}}>\veps \right)\\
&\le \Pp\left(n^{-1/3}\sum_{i=1}^{n} \wt{i}\I{J_i \le tn^{2/3}, \hat{E}_i \le 2sn^{1/3}}>\veps \right)\\
&\quad+ \Pp(\hat{E}_{\sroot_{\ast}}>2sn^{1/3}),  
\end{align*} 
where the first inequality follows from the fact that the first $\lfloor sn^{1/3}\rfloor$ roots are necessarily among those $i$ with $\hat{E}_i\le \hat{E}_{\sroot_{\ast}}$.
We see that the first term goes to $0$ as $n\to \infty$ by Lemma~\ref{lem:competing_clocks}. For the second term, recall the definition of the event $E(k,r,x,y,z)$ from the proof of Proposition~\ref{prop:roots_uar}, and in particular that it implies that $\hat{E}_{\sroot_k}\le z$. Moreover, recall that for $k=\lfloor sn^{1/3}\rfloor$ and $z=2sn^{1/3}$ we may pick $r,x,z$ such that $\Pp(E(k,r,x,y,z)^c)<\veps$ for $n$ large enough, and the second statement of the proposition follows. The first statement follows completely analogously. 

For the fourth statement, observe that 
\begin{align*} 
&\Pp\left(n^{-1/3}\sum_{i=1}^{\lfloor s n^{1/3}  \rfloor} \I{J_{\hat{R}_i} \le tn^{2/3}}>\veps \right)\\
&\qquad \le \Pp\left(n^{-1/3}\sum_{i\in[n]} \I{J_i \le tn^{2/3}, \hat{E}_i \le 2sn^{1/3}}>\veps \right)+ \Pp(\hat{E}_{\hat{R}_{\ast}}>2sn^{1/3})\end{align*}
The first term goes to $0$ as $n\to \infty$ by Lemma~\ref{lem:competing_clocks} and we controlled the second term above. The third statement of the proposition follows analogously. 
\end{proof}

Recall that $\wc(s) = \sup\{i:\wst_i \leq s\}$ is the index of the component that is being explored at time $s$ when we select roots in a weight-biased manner, while $\scc(s) = \sup\{i:\sst_i \leq s\}$ is the index of the component that is being explored at time $s$ when we select roots uniformly. The following lemma is the last ingredient we need for the proof of~\cref{prop: comps explored}.

\begin{lemma}\label{lem:analysis_hitting_time}
Suppose $f_n(t,s)$ is a function from $\R_+ \times \R_+$ to $\R_+$ that is non-decreasing in $t$ for all $s$, that $f(t,s)$ is a function from $\R_+ \times \R_+$ to $\R_+$ that is constant in $s$, non-decreasing and continuous in $t$ and that $f_n\to f$ uniformly on compact sets. Moreover, suppose that $g_n(s)$ is a function from $\R_+ $ to $\R_+$ that is non-decreasing, that $g(s)$ is a function from $\R_+ $ to $\R_+$ that is strictly increasing, continuous, satisfies $g(0)=0$ and $\lim_{t\to\infty}g(s) = \infty$, and that $g_n\to g$ uniformly on compact sets. Then, 
\[\left(\sup\{s:g_n(s) \le f_n(t,s)\}, t\ge 0\right) \to \left(\sup\{s:g(s) \le f(t,s)\}, t\ge 0\right)\]
uniformly on compact sets. 
\end{lemma}
\begin{proof}
Since $f(t,s)$ is constant in $s$, we write $f(t)$ to simplify notation in the proof, and define the functions
\[
h(t) := \sup\{s:g(s) \le f(t)\} \qquad \text{and} \qquad h_n(t):= \sup\{s:g_n(s) \le f_n(t,s)\}, \qquad t \ge 0.
\]
By our assumptions on $g$, it is a bijection from $\R_+$ to $\R_+$, and therefore invertible. It follows that $h(t)=g^{-1}(f(t))$ for all $t\ge 0$. Let $t\ge 0$ and we will show that $h_n(t)\to h(t)$. 
Let us first argue that $\liminf_{n\to\infty} h_n(t)\ge h(t)$.
We note that if $f(t)=0$, then $h(t)=0$, so that $h_n(t)\ge h(t)$ holds automatically. 
We can therefore assume $f(t)>0$, so that $h(t)>0$. Let us take $0\le s<h(t)$. Then we have $g(s)<f(t)$ by the monotonicity of $g$.  
Take $\veps' = (f(t) - g(s))/3. $ Since $g_n \rightarrow g$ and $f_n \rightarrow f$ uniformly on compact sets, there exists $N$ large enough such that $\left\vert g_n(s) - g(s) \right\vert < \varepsilon'$ and $\left\vert f_n(t,s) - f(t) \right\vert < \varepsilon'$ for all $n \geq N$. It follows that 
\[
g_n(s) < g(s) + \varepsilon' < f(t) - \varepsilon' < f_n(t,s),
\]
which implies $s\le h_n(t)$ for all $n \geq N$. This proves $\liminf_{n\to\infty}h_n(t)\ge s$. Letting $s$ increase to $h(t)$ yields $\liminf_{n\to\infty}h_n(t)\ge h(t)$. We have now completed the proof that $\liminf_{n\to\infty}h_n(t)\ge h(t)$. To show that $\limsup_{n\to\infty}h_n(t)\le h(t)$, we can argue in a similar way by taking $s>h(t)$. We have shown that $h_n(t)\to h(t)$ pointwise. The uniform convergence readily follows, since both $h_n$ and $h$ are monotone and $h$ is further continuous by assumption.
\end{proof}

\begin{proof}[Proof of~\cref{prop: comps explored}]
We will show the convergence of $\wc$ under rescaling. The proof for the convergence of $\scc$ is analogous. 

Fix $T,U>0$. From~\eqref{eq:tauk}, we see that $\tau_j\le t$ if and only if $\inf_{0\le u\le  t}S_n(u)\le -\sum_{i\in [j]}\wt{\wroot_i}\I{J_{\wroot_i}> t}$. We then deduce that
\[\wc(t)=\sup\left\{j: \sum_{i=1}^j \wt{\wroot_i} \le -\inf_{0\le u \le t} \left( S_n(u) -\sum_{i=1}^j \wt{\wroot_i}\I{J_{\wroot_i} \le u}\right) \right\}. \]

By Proposition~\ref{prop:roots_uar},
\[\left(n^{-1/3}\sum_{i=1}^{\lfloor s n^{1/3} \rfloor} \wt{\wroot_i}, s\ge 0\right) \convprobn (s, s\ge 0)\]
uniformly on compact sets. Moreover, by the convergence of $S_n$ to $\bm$ under rescaling in Proposition~\ref{prop:conv_S_n}, and Proposition~\ref{prop:corr_uar}, we see that 
\begin{align*} &\left(n^{-1/3}\inf_{0\le u \le tn^{2/3}} \left( S_n(u)-\sum_{i=1}^{\lfloor s n^{1/3} \rfloor} \wt{\wroot_i}\I{J_{\wroot_i} \le u}\right), 0\le t\le T, 0\le s \le U \right) \\
&\qquad \convdist \left( \inf_{0 \le u \le t} W_u,  0\le t \le T, 0\le s \le U \right) \end{align*}
uniformly, jointly with the convergence of $S_n$ to $\bm$ under rescaling in Proposition~\ref{prop:conv_S_n}. 
This allows us to apply Lemma~\ref{lem:analysis_hitting_time} with
\[
g_n(s) = n^{-1/3}\sum_{i=1}^{\lfloor s n^{1/3} \rfloor} \wt{\wroot_i}, \quad f_n(t, s) = -n^{-1/3}\inf_{0\le u\le tn^{2/3}}\left( S_n(u)-\sum_{i=1}^{\lfloor s n^{1/3} \rfloor} \wt{\wroot_i}\I{J_{\wroot_i} \le u}\right),
\]
as well as $g(s) =s$, $f(t, s)=-\inf_{0\le u\le t}W_u$. We note that in this case,
\[
n^{1/3}\sup\{s: g_n(s)\le f_n(t, s)\}= c(tn^{2/3})+1,
\]
and for each $t\ge 0$, 
\[
\sup\{s: g(s)\le f(t, s)\} = -\inf_{0\le u\le t}W_u. 
\]
The conclusion follows.
\end{proof}

\section{Proof of \cref{lem:eps_deldisguise}}\label{proof of disguised delta eps}
Rather than proving \cref{lem:eps_deldisguise} directly, we will prove and apply a more general result, for the statement of which we require the following definitions. Recall that for $f \in \Dt$, we let $\exc{i}(f)=(g^{(i)}(f),d^{(i)}(f))$ be the $i$th longest excursion of $f$. 
\begin{definition}[good function]\label{good function}
For $T>0$, a continuous function $f:[0,T]\to \R$ with $f(0)=0$ is {\em good} if the following conditions are satisfied.
\begin{enumerate}[(G1)]
    
    \item\label{cond 4} For any excursion $(g,d)$ of $f$, for any $\delta>0$, 
    \[
    \inf_{s< d+\delta}f(s) < f(d);
    \]
    {\item\label{all in excursions} 
    Let $T^-=\sup\{t \in [0,T]: f(t)=\underline{f}(T)\}$; then $[0,T^-]\setminus \bigcup_{i\ge 1} [g^{(i)}(f),d^{(i)}(f)]$ has Lebesgue measure zero;}
    
    \item\label{fix k'} The ordered excursion lengths $(\excsize{i}(f),i \ge 1)$ are strictly decreasing.
\end{enumerate}
\end{definition}

For a function $h:[0,T] \to \R$ and $\veps>0$, let  $C(h,\veps)$ be the set of all $\veps$-dense sets of cutoff levels $\alpha$ such that 
\begin{equation}
\label{cond: cutoff}
- \max\{\alpha_i:\alpha_i\in \alpha\} < \underline{h}(T)+ \veps.
\end{equation}
For $\veps>0, k\in \N$, a good function $f$ and
$h \in \mathbb{D}_+[0,T]$, we define
\begin{equation*}
\Delta_\veps(f,h; k) \coloneqq \sup_{\alpha \in C(h,\veps)} \sup_{i\in[k]} \left(|g^{(i,\alpha)}(h) - g^{(i)}(f)| + |d^{(i,\alpha)}(h) - d^{(i)}(f)|\right),
\end{equation*}
where we recall that $g^{(i,\alpha)}(h)$ and $d^{(i,\alpha)}(h)$ denote the left and right endpoints, respectively, of the $i$th largest $(h,\alpha)$-excursion.

\begin{lemma}\label{delta epsilon}
Let $\delta>0, k\geq 1$ and $T>0$. Then for any good function $f :[0,T]\to \R$ such that that $0<g^{(i)}(f)<d^{(i)}(f)<T$ for all $i\in [k]$, there is $\veps=\veps(\delta,T,k,f)>0$ such that for any $h \in \Dt$ with $h(0)=0$ and $\|f - h\|_\infty < \veps$, 
it holds that $\Delta_\veps(f,h;k) < \delta$. 
\end{lemma}

\cref{lem:eps_deldisguise} follows almost immediately from \cref{delta epsilon}.

\begin{proof}[Proof of \cref{lem:eps_deldisguise}]
Let $(B(t),t \in [0,T])$ be a standard Brownian motion. Then the law of $(W(t),t \in [0,T])$ is absolutely continuous with respect to that of $(B(t),t \in [0,T])$, and it follows from standard properties of Brownian motion that $(W(t),t\in[0,T])$ is almost surely a good function, and that almost surely $0<g^{(i)}<d^{(i)}<T$ for all $i\in [k]$. 

Next, using Skorohod's representation theorem, we may 
assume that 
\[\mathbf{P}\big(\|\overline{X}_n-W|_{[0,T]}\|_\infty\to 0\big)=1,\]
or equivalently, that
\[\mathbf{P}\big(\|\overline{X}_n-W|_{[0,T]}\|_\infty<\varepsilon,\forall\varepsilon>0,\forall\,n\mbox{ sufficiently large}\big)=1.
\] 
Thus, given $k\ge 1$ and $\delta>0$, by \cref{delta epsilon} there exists $\varepsilon=\varepsilon(\delta,k)>0$ such that 
\[
\begin{aligned}1&=\mathbf{P}\big(\|\overline{X}_n-W|_{[0,T]}\|_\infty<\varepsilon,\forall\,n\mbox{ sufficiently large}\big)
\\
&\le\mathbf{P}\big(\Delta_{\veps}(W|_{[0,T]},\overline{X}_n;k)<\delta,\forall\,n\mbox{ sufficiently large}\big).\end{aligned}
\]

For $i\ge 1$ and a set of cutoff levels $\alpha$, let $(g^{(i,\alpha)},d^{(i,\alpha)}) := \exc{i,\alpha}(\overline{X}_n)$ and $(g_T^{(i)},d_T^{(i)}) \coloneqq \exc{i}(\bm|_{[0,T]})$.
With $\alpha^{(n)}$ and $\beta^{(n)}$ fixed as above, it follows that for any $\delta>0$ and $k\ge 1$, for all $n$ large enough,
\begin{align*}
&\max_{i \in [k]} 
\left(
|g^{(i,\alpha^{(n)})}-g^{(i)}_T|+
|d^{(i,\alpha^{(n)})}-d^{(i)}_T|+
|g^{(i,\beta^{(n)})}-g^{(i)}_T|+
|d^{(i,\beta^{(n)})}-d^{(i)}_T|
\right)
\\
&\le 2\Delta_\varepsilon\big(W|_{[0,T]},\overline{X}_n;k\big)
<2\delta.
\end{align*}
This gives the desired joint convergence in the product topology.
\end{proof}

\begin{proof}[Proof of \cref{delta epsilon}]
Fix $k\geq 1$ and $T>0$. Let $f:[0,T]\to\mathbb{R}$ be a good function. 
Then there exists $\delta_0>0$ such that for all $2\le i\le k$, 
\begin{equation}\label{exc f lower bd}
\excsize{i}(f) < \excsize{i-1}(f) - \delta_0.
\end{equation}

We let $(a_1,b_1),\ldots, (a_{m'},b_{m'})$ denote the intervals of $[0,T^-] \backslash \bigcup_{i\in[k]}E^{(i)}(f)$. We let $\intv_1,\ldots, \intv_{m}$ denote the intervals of the set 
\[
\bigcup_{i\le m'}
\Big(a_i - \frac{\delta_0}{4}, b_i + \frac{\delta_0}{4}\Big)\cap [0,T^-],
\]
noting that $m\le m'$. 
Observe that for $j\in [m]$, every 
$f$-excursion contained in $\intv_j$ has length less than $\excsize{k}(f)- \delta_0$ by \labelcref{exc f lower bd}. Hence, since $f$ satisfies \labelcref{all in excursions}, we can find a set of cutoff levels 
$\beta_j=\{\beta_{1,j}<\cdots<\beta_{l_j,j}\} \subset \intv_j$ 
which is $\frac12(\excsize{k}(f)- \delta_0)$-dense in $\intv_j$ such that $\beta_{i,j}$ is an endpoint of an $f$-excursion for all $i\in[l_j]$. 

We abbreviate $g_i=g^{(i)}(f)$ and $d_i=d^{(i)}(f)$ for $i\in [k]$.
Fix $\delta>0$, and let {$\delta'>0$ be small enough so that $\delta'\le\min\{\frac{\delta}{2},\frac{\delta_0}{4}\}$ and $\delta'<g_i<g_i+\delta'<d_i-\delta'<d_i<T-\delta'$ for all $i\in [k]$}. By \labelcref{cond 4} and continuity of $f$, 
there exists $\veps=\varepsilon(\delta')>0$ such that for all $i\in[k]$,
\begin{equation}
\begin{aligned}
   \inf_{s\leq g_i - \delta'}f(s) > f(g_i) + 5\veps,\hspace{.5cm}\inf_{s\leq d_i+\delta'}f(s) < f(d_i) - 4\veps, 
\label{cond 4 in proof}
\end{aligned}
\end{equation}
and 
\begin{equation}
\begin{aligned}
   \inf_{s\in(g_i+\delta',d_i-\delta')}f(s) &> f(d_i) + 4\veps,
\label{lower bd in proof}
\end{aligned}
\end{equation}
and such that for all $d\in\{\beta_{i,j}: i\in [\ell_j], j\in [m]\}$, 
\begin{equation}\label{f continues decreasing}
    \inf_{s\in (d,  d+\delta')}f(s) < f(d)-4\veps.
\end{equation}
Let $h \in \mathbb{D}_+[0,T]$ be such that $h(0)=0$ and 
\begin{equation}\label{f h similar}
    \|f-h\|_\infty<\veps,
\end{equation}
and let $\alpha \in C(h,\veps)$ be an $\veps$-dense sequence of cutoff times. We proceed in steps. 

\textit{Step~1. } Recall from (\ref{tauhz}) that $\tau(h,\alpha_i)=\inf\{t \ge 0: h(t)=-\alpha_i\}\wedge T$. We first show that for all $j\in [k]$, there is an 
$i_j$ such that 
\begin{equation}
\label{exc approximates}
   \big|\tau(h,\alpha_{i_j-1})- g_j| < \delta'\quad\text{and}\quad |\tau(h,\alpha_{i_j})-d_j\big| < \delta'.
\end{equation}
We prove \labelcref{exc approximates} in the case $j=1$; the argument works identically for other $j \in [k]$.
{
Let us first argue that
\begin{equation}
\label{def: Adel}
A(\delta'):=\{\alpha'\in \alpha: |\tau(h, \alpha')-g_1|<\delta'\} \ne \varnothing \quad\text{and}\quad \max_{\alpha'\in \alpha}\alpha'\notin A(\delta'). 
\end{equation}
To prove~\eqref{def: Adel}, we note that 
\begin{equation}
\label{bd: h}
h(g_1)>f(g_1)-\veps = f(d_1)-\veps\ge \inf_{s\le d_1+\delta'}f(s)+3\veps\ge \underline{f}(T)+3\veps\ge \underline{h}(T)+2\veps,
\end{equation}
by~\eqref{f h similar} and~\eqref{cond 4 in proof}. Since $\alpha$ is $\veps$-dense and, by~\eqref{cond: cutoff}, $- \max\{\alpha_i:\alpha_i\in \alpha\} < \underline{h}(T)+ \veps$, 
there exists some $\alpha_{\ast}\in \alpha$ such that 
\begin{equation}
\label{alpha in interval}
   h(g_1) \le -\alpha_{\ast} \le h(g_1) + 2\veps.
\end{equation} 
We apply~\eqref{f h similar} and~\eqref{cond 4 in proof} again to find that
\[
\inf_{s\le g_1-\delta'}h(s)\ge \inf_{s\le g_1-\delta'}f(s)-\veps \ge f(g_1)+4\veps\ge h(g_1)+3\veps,
\]
Combined with~\eqref{alpha in interval}, this yields that $\inf_{s\le g_1-\delta'}h(s)>-\alpha_{\ast}\ge h(g_1)$. Since $h$ has no negative jumps, we have $\tau(h, \alpha_{\ast})\in (g_1-\delta, g_1]$. This shows the set $A(\delta')$ is non empty. To see why $\max\{\alpha':\alpha\in\alpha\}\notin A(\delta')$, we again use that $-\max\{\alpha':\alpha\in\alpha\}\le \underline{h}(T)+\veps$. However, since $g_1+\delta'<d_1$, we have $\inf_{s\le g_1+\delta'}f(s)\ge f(d_1)> \underline{h}(T)+3\veps$ according to~\eqref{bd: h}; so if $\alpha'\in A(\delta')$, then we would have $-\alpha'\ge \underline{h}(T)+3\veps$.
This proves~\eqref{def: Adel}.
Thanks to it, we see that  $\max\{\alpha': \alpha'\in A(\delta')\}$ cannot be the last cutoff level. So we can choose $i_1$ such that $\alpha_{i_1-1}=\max\{\alpha': \alpha'\in A(\delta)\}$. }

By the definition of $\alpha_{i_1-1}$, we now have $\tau(h,\alpha_{i_1}) > g_1+\delta'$. 
If $\tau(h,\alpha_{i_1}) \leq d_1 - \delta'$, then $f(\tau(h,\alpha_{i_1})) > f(g_1) + 4\veps$ by \labelcref{lower bd in proof}, so since $\|f-h\|_\infty < \veps$ it follows that 
\begin{equation}
\label{bad alpha l}
   -\alpha_{i_1} = h\big(\tau(h,\alpha_{i_1})\big)  > h(g_1) + 2\veps.
\end{equation}
However, $\alpha_{i_1-1}$ and $\alpha_{i_1}$ were chosen so that $-\alpha_{i_1} < -\alpha_{i_1-1} \le -\alpha_{\ast}<h(g_1) +2\veps$, 
contradicting \labelcref{bad alpha l}.
Hence $\tau(h,\alpha_{i_1}) > d_1-\delta'$. 

It remains to show that $\tau(h,\alpha_{i_1}) < d_1 + \delta'$. By \eqref{cond 4 in proof} and \eqref{f h similar}, we have $h(d_1)>f(d_1) - \veps$ and $h(d_1+\delta')<f(d_1)-3\veps$. Thus, $h(d_1) - h(d_1 + \delta')>2\veps$, which means that there exists a cutoff level $\alpha'\in \alpha$ such that $\alpha'\in (h(d_1 + \delta'),h(d_1))$ and $h$ hits it in the interval $(d_1,d_1+\delta')$. This proves that $\alpha_{i_1}<d_1+\delta'$.

\textit{Step~2. } {Recalling that $\tau(h, \alpha_i)-\tau(h,\alpha_{i-1})$ corresponds to the length of an $(h, \alpha)$-excursion,} we next show that for all $j\in[k-1]$, 
\begin{equation}\label{h excs are ordered}
   \tau(h,\alpha_{i_j})-\tau(h,\alpha_{i_{j}-1})>\tau(h,\alpha_{i_{j+1}})-\tau(h,\alpha_{i_{j+1}-1}).
\end{equation}

For all $1 \le j <j'\leq k$, it holds that
\begin{align*}
\tau(h,\alpha_{i_{j}})-\tau(h,\alpha_{i_{j}-1})&\geq \excsize{j}(f)-2\delta'\quad\text{by ~\eqref{exc approximates}}\\ 
&> \excsize{j'}(f) -2\delta' + \delta_0 \quad\text{by \labelcref{exc f lower bd}}\\
&>  \tau(h,\alpha_{i_{j'}})-\tau(h,\alpha_{i_{j'}-1}) - 4\delta'+\delta_0 \quad\text{by ~\eqref{exc approximates}}\\
&\ge \tau(h,\alpha_{i_{j'}})-\tau(h,\alpha_{i_{j'}-1})
\quad \text{since $\delta_0 \ge 4\delta'$},
\end{align*}
so the excursions are indeed ordered by size. 

\textit{Step 3. } Recall that $\exc{j,\alpha}(h)$ denotes the $j$th largest of the excursions $(h, \alpha)$-excursions.
We finally show that for all $j\in[k]$,
\begin{equation}
\label{Step3}
   \exc{j, \alpha}(h) = \big(\tau\big(h,\alpha_{i_{j}-1}), \tau\big(h,\alpha_{i_{j}})\big).
\end{equation}
Suppose towards a contradiction that this is not the case, and let $j\coloneqq \min\{j' \in [k]:\, \excsize{j',\alpha}(h) > 
\tau(h,\alpha_{i_{j'}})-\tau(h,\alpha_{i_{j'}-1})\}$. Then $\excsize{j,\alpha}(h)\neq \tau(h,\alpha_{i_{j'}})-\tau(h,\alpha_{i_{j'}-1})$ for {all} $j'\in[k]$, {since we have seen previously that the excursions lengths are strictly decreasing. Combining this with~\eqref{exc approximates}, we deduce that} $\badexc$ intersects any of the $f$-excursions $\exc{1}(f),\ldots,\exc{k}(f)$ by at most $\delta'\le \frac{\delta_0}{4}$ and can intersect at most two such excursions. It follows that
\begin{equation}
\label{small exc overlap}
\big|\badexc(h) \cap \bigcup_{i=1}^k \exc{i}(f)\big|\leq 2\delta' \le \frac{\delta_0}{2}.
\end{equation}
Furthermore, 
\begin{align*}
\badexcsize(h) &\ge \tau\big(h,\alpha_{i_{k}}\big)-\tau\big(h,\alpha_{i_{k}-1}\big)
\\
&\ge \excsize{k}(f) - 2\delta' \quad\text{by \labelcref{exc approximates}}\\
&\ge \excsize{k}(f) - \frac{\delta_0}{2}. \stepcounter{equation}\tag{\theequation}\label{bad exc lower bound}
\end{align*}
Recall the intervals $\intv_1,\ldots, \intv_m$ defined earlier in the proof. It follows from \labelcref{small exc overlap} and \labelcref{bad exc lower bound} that $\badexc(h)\subset \intv_{j'}$ for some $j'\leq m$. Recall that for all $j\le m$, $\beta_j=\{\beta_{1,j},\ldots,\beta_{l_{j},j}\}$ is a set of cutoff levels that is $\frac12(\excsize{k}(f)-\delta_0)$-dense in $\intv_{j}$ by construction. To simplify notation below, we define $(g,d) \coloneqq (g^{(j,\alpha)}(h), d^{(j,\alpha)}(h))$.
Then, by \eqref{bad exc lower bound},
\[
|(g,d-\delta')| = \badexcsize(h) - \delta'   >  \excsize{k}(f) - \delta_0,
\]
so there exists $i \in [\ell_j]$ such that $\beta_{i,j'}\in (g,d-\delta')$. Hence, by \labelcref{f continues decreasing}, there is $s \in (d-\delta',d)$ such that $
f(s)<f(\beta_{i,j'}) - 4\veps$, so $h(s) < f(\beta_{i,j'}) - 3\veps$ by \labelcref{f h similar}. Therefore, there exists a cutoff level $\alpha'_{\ast} \in \alpha$ 
such that 
\begin{equation}\label{bound 10}
h(s) < -\alpha'_{\ast} < h(s) +2\veps  < f(\beta_{i,j'}) - \veps.
\end{equation}
Since $\beta_{i,j'}$ is the endpoint of an $f$-excursion, it holds that $f(g) \ge f(\beta_{i,j'})$, so $h(g) \ge f(\beta_{i,j'}) - \veps > -\alpha'_{\ast}$ by \labelcref{f h similar} and \labelcref{bound 10}. Since $h$ does not have any negative jumps, and $h(g)$ is itself a cutoff level for $h$, this implies that $\tau(h,\alpha'_{\ast}) \in(g,d)$. Thus, we have shown that $(g,d)$ is not of the form $(\tau(h,\alpha_{i'-1}),\tau(h,\alpha_{i'}))$ for some $i'\in[k]$. This yields a contradiction and proves (\ref{Step3}). 

Combining \eqref{exc approximates},\eqref{h excs are ordered} and \eqref{Step3} gives that
\begin{equation*}
\sup_{i\in [k]}\Big(\big|{g^{(i,\alpha)}(h)} - g_i\big| + \big|{d^{(i,\alpha)}(h)} - d_i\big|\Big) < 2\delta'.
\end{equation*}
Since $\alpha \in C(h, \veps)$ was arbitrary, we can take the supremum over all such sets of cutoff levels to conclude that
\begin{equation*}
\Delta_{\veps}(f,h;k)=\sup_{\alpha \in C(h, \veps)} \sup_{i\in [k]} \Big(\big|{g^{(i,\alpha)}(h)} - g_i\big| 
+ \big|{d^{(i,\alpha)}(h)} - d_i\big|\Big) < \delta.\qedhere
\end{equation*}
\end{proof}

\appendix

\section{Technical results}
\label{sec: comment}

\begin{proof}[Proof of~\cref{prop:conv_S_n}]
Denoting by $q_i$ the rate of $J_i$, we have the following Doob--Meyer decomposition of $S_n$: 
\[
S_n = A_n + M_n, \text{ where } A_n(t)=-t+\sum_{i\in[n]}w^{(n)}_i q_i\int_0^t\I{J_i>s}ds,
\]
and $M_n$ is a zero-mean martingale with $\mathbb E[M_n^2(t)]=\sum_i (w^{(n)}_i)^2\Pp(J_i\le t) $. 
Recall from~\eqref{defW} the definition of $W$. We aim to show that for any $T>0$, 
\begin{equation}
\label{cv: AnMn}
\Big(n^{-1/3}A_n(n^{2/3}t), n^{-1/3}M_n(n^{2/3}t), t\ge 0\Big) \convdistn \Big(\lambda t-\frac{\mu't^2}{2\mu^2}, \sqrt{\frac{\mu'}{\mu}}B(t), t\ge 0\Big)
\end{equation}
uniformly on $[0, T]$, 
which will then imply the convergence of $\scaledrw$. To that end, we note that $A_n$ is differentiable in $t$, and that for $t \ge 0$ we have
\[
n^{-1/3}A_n(n^{2/3}t)= \int_0^t a_n(s)ds \quad\text{with}\quad a_n(s)=-n^{1/3}+n^{1/3}\sum_{i\in [n]}\wt{i}q_i\I{J_i>n^{2/3}s}.
\]
We note that $1-\exp(-x)=x+r(x)\cdot x^2$, where the function $r$ is uniformly bounded on any compact set of $\R_+$. Recalling that $q_i = \wt{i}(1+\cw n^{-1/3})/|\wtv|_1$, we also have 
    \begin{equation}
    \label{bd: r}
    \sup_{n\ge 1}\max_{i\in [n]}q_i n^{2/3} t\le (1+|\lambda|) t\, \frac{n^{2/3}\cdot\max_{i\in [n]}\wt{i}}{|\wtv|_1}<\infty,
    \end{equation}
    thanks to \labelcref{maximal weight} and \labelcref{first moment}. Moreover, 
    \begin{align*}
    \mathbb E[a_n(t)] &= -n^{1/3} + n^{1/3}\sum_{i\in [n]}\wt{i}q_i\Pp(J_i>n^{2/3}t) \\
    &= - n^{1/3} + n^{1/3}\sum_{i\in [n]} \wt{i}q_i e^{-q_i n^{2/3}t} \\
    & = n^{1/3}\Big(\sum_{i\in [n]}\wt{i}q_i-1\Big)-n^{1/3}\sum_{i\in [n]} \wt{i}q_i (1-e^{-q_i n^{2/3}t})\\
    & = n^{1/3}\Big(\sum_{i\in [n]}\wt{i}q_i-1\Big) - nt\sum_{i\in [n]}\wt{i}q_i^2-n^{5/3}t^2\sum_{i\in [n]}\wt{i}q_i^3r\big(q_in^{2/3}t\big)
    \end{align*}
Thanks to~\eqref{bd: r}, we can find some $C\in (0,\infty)$ so that 
\begin{align}\notag
n^{5/3}t^2\sum_{i\in [n]}\wt{i}q_i^3 \,r\big(q_in^{2/3}t\big)&\le Cn^{5/3}t^2\sum_{i\in [n]}\wt{i}q_i^3\\ \label{cv: error}
&\le C(1+|\lambda|)^3t^2 n^{5/3} \frac{\max_i \wt{i}\cdot|\wtv|_3}{|\wtv|_1^3}\xrightarrow{n\to\infty} 0,
\end{align}
due to the assumptions~\eqref{maximal weight},~\eqref{first moment} and~\eqref{third moment}. We further deduce, with the additional assumption~\eqref{second moment}, that
\[
n^{1/3}\Big(\sum_{i\in [n]} \wt{i}q_i-1\Big)\xrightarrow{n\to\infty} \lambda  \quad\text{and} \quad nt\sum_{i\in [n]}\wt{i}q_i^2\xrightarrow{n\to\infty}\frac{\mu'}{\mu^2 }t.
\]
Combined with~\eqref{cv: error}, this shows that  
\[
\mathbb E[a_n(t)] \xrightarrow{n\to\infty} \lambda -\frac{\mu'}{\mu^2}t. 
\]
On the other hand, using the independence of the collection $J_i, i\in [n]$ and the bound $1-e^{-x}\le x$ for $x\ge 0$, we find that
\begin{align*}
\mathrm{Var}\big(a_n (t)\big)&\le n^{2/3}\sum_{i\in [n]}(\wt{i})^2q_i^2\,\Pp(J_i\le n^{2/3}t)\le n^{4/3}t\sum_{i\in [n]}(\wt{i})^2q_i^3\\
&\le (1+|\lambda|)^3 t\,n^{4/3}\frac{(\max_i\wt{i})^2\cdot|\wtv|_3}{|\wtv|_1^3}\xrightarrow{n\to\infty}0, 
\end{align*}
where the final convergence follows from~\eqref{maximal weight},~\eqref{first moment} and~\eqref{third moment}. Together with the previous convergence of its expectation, this shows that for each $t\ge 0$, 
\begin{equation}
\label{cv: an-t}
a_n(t) \convprobn \lambda -\frac{\mu'}{\mu^2}t.
\end{equation}
Using the monotonicity of $a_n$, let us now explain how the convergence in~\eqref{cv: an-t} can hold uniformly on $[0, T]$. We first partition $[0,T]$ into intervals of length $\delta$, get a uniform bound on each of these, then patch them together. Let $\delta\in (0, 1)$ be such that $(1+\mu'/\mu^2)\delta<\veps$, and let 
$t_i = i\delta \wedge T$ for $i\ge 0$. 
We deduce from~\eqref{cv: an-t} that for $n$ large enough and $1\le i\le \lfloor T/\delta\rfloor+1$, 
    \begin{equation}\label{eq:tube''}
    \Pp\left(\Big|a_n(t_i)- \lambda+\frac{\mu'}{\mu^2}t_i\Big|>\delta\right) < \frac{\veps}{2T}.
    \end{equation}
    For any $i\ge 0$ and $t\in [t_i,t_{i+1}]$, by monotonicity of $a_n(t)$ in $t$,
    \[
    a_n(t_{i+1})-\lambda + \frac{\mu'}{\mu^2}t_{i+1}-\frac{\mu'}{\mu^2}\delta \le a_n(t)-\lambda+\frac{\mu'}{\mu^2}t\le a_n(t_{i})-\lambda + \frac{\mu'}{\mu^2}t_{i}+\frac{\mu'}{\mu^2}\delta ,
    \]
    and hence, using $\eqref{eq:tube''}$ and a union bound over $1\le i \le 2T/\veps$, 
    \[
    \Pp\Big(\sup_{t\in[0, T]}\Big|a_n(t) - \lambda+ \frac{\mu'}{\mu^2}t\Big| > \veps\Big)<\veps, 
    \]
    for $n$ large enough. 
Integrating over $t$, we then obtain the uniform convergence of $A_n$ as claimed in~\eqref{cv: AnMn}. 

For the convergence of $M_n$, we use the martingale central limit theorem. Indeed, it holds that $\mathbb E[M_n(t)]=0$ and 
\begin{align*}
n^{-2/3}\mathbb E\big[M_n^2(n^{2/3}t)\big] & = n^{-2/3}\sum_{i\in [n]}(\wt{i})^2\Pp(J_i\le n^{2/3}t) = n^{-2/3}\sum_{i\in[n]}(\wt{i})^2(1-e^{-q_in^{2/3}t}) \\
&= t\sum_{i}(\wt{i})^2q_i + n^{2/3}t^2\sum_{i\in [n]}(\wt{i})^2q_i^2\, r(q_i n^{2/3}t).
\end{align*}
With a similar argument as previously, we deduce from the assumptions of $\wtv$ and~\eqref{bd: r} that 
\[
n^{-2/3}\mathbb E\big[M_n^2(n^{2/3}t)\big] \to \frac{\mu'}{\mu}t. 
\]
Moreover, since $A_n$ is continuous in $t$, we have
\[
\sup_{t\in[0, Tn^{2/3}]}n^{-1/3}\big|M_n(t)-M_n(t-)\big| = \sup_{t\in[0, Tn^{2/3}]}n^{-1/3}\big|S_n(t)-S_n(t-)\big|  \le n^{-1/3}\max_{i\in[n]}\wt{i} \to 0,
\]
by~\eqref{maximal weight}. According to~\cite[Chapter 7, Theorem 1.4]{Markov-book}, the convergence in~\eqref{cv: AnMn} follows. 
\end{proof}
For the proof of~\cref{LLN lemma} as well as~\cref{lem:lln_roots}, we introduce the following process. For each $n\in \N$, let $\mathbf x_n=(x_n(i))_{i\in [n]}$ and $\mathbf y_n=(y_n(i))_{i\in [n]}$ be two sequences of positive real numbers, and denote
\[
Q_n(t) = \sum_{i\in [n]}x_n(i)\I{\mathcal E_i\le t}, \quad t\ge 0, 
\]
where $(\mathcal E_i, i\in [n])$ is a collection of independent exponential variables of respective rates $y_n(i)$. 
\begin{lemma}
\label{lem: LLN-gen}
Suppose that there exists $\alpha\in (0, \infty)$ such that
\begin{equation}
\label{cond: LLN-gen}
\sum_{i\in [n]}x_n(i)y_n(i)\xrightarrow{n\to\infty} \alpha \quad\text{and}\quad n^{2/3}\max_{i\in [n]}y_n(i) \xrightarrow{n\to\infty} 0. 
\end{equation}
Suppose further that $(b_n)_{n\in \N}$ is a sequence of positive numbers that satisfies  $b_n=O(n^{2/3}), n\to\infty$ and
\begin{equation}
\label{cond: bn}
b_n^{-1}\max_{i\in [n]}x_n(i) \xrightarrow{n\to\infty}0.
\end{equation}
Then we have for any $\veps, T>0$, 
\[
\Pp\Big(\sup_{t\in [0, T]}\big|b_n^{-1}Q_n(b_n t)- \alpha t\big|>\veps\Big)\xrightarrow{n\to\infty} 0.  
\]
\end{lemma}

\begin{proof}
As in the previous proof, we note that $1-\exp(-x)=x+r(x)\cdot x^2$, where the function $r$ is uniformly bounded on any compact set of $\R_+$. We have 
\begin{align*}
        b_n^{-1}\mathbb{E}[Q_n(b_nt)] &= b_n^{-1}\sum_{i\in[n]}x_n(i)\Pp(\mathcal E_i \le b_nt)\\ 
        &= t\sum_{i\in [n]}x_n(i)y_n(i) + b_nt\sum_{i\in [n]}x_n(i)y^2_n(i)r\big(y_n(i) b_n t\big) \\
        & = \alpha t(1+o(1))+b_nt\max_{i\in [n]}y_n(i) \sum_{i\in [n]}x_n(i)y_n(i)r\big(y_n(i) b_n t\big)\\
        &\xrightarrow{n\to\infty}\alpha t,
\end{align*}
thanks to~\eqref{cond: LLN-gen}, the assumption that $b_n=O(n^{2/3}t)$ and the fact that $\sup_{n\ge 1}\sup_{i\in[n]}r(y_n(i)b_nt)<\infty$. On the other hand, using the independence of the collection $(\mathcal E_i, i\in [n])$, the bound $1-e^{-x}\le x$ for $x\ge 0$, we find that
    \begin{align*}
        b_n^{-2}\mathrm{Var}\big(Q_n(b_n t)\big) &=  b_n^{-2}\sum_{i\in [n]} x_n^2(i)\mathrm{Var}(\I{\mathcal E_i\le b_nt}) 
        \le b_n^{-2}\sum_{i\in [n]}x_n^2(i)\Pp(\mathcal E_i\le b_nt) \\
        & \le b_n^{-1}t\sum_{i\in [n]}x_n^2(i)y_n(i) \le t b_n^{-1}\max_{i\in [n]}x_n(i)\sum_{i\in [n]}x_n(i)y_n(i)\to 0,
    \end{align*}
thanks to~\eqref{cond: LLN-gen} and~\eqref{cond: bn}. It then follows from Chebyshev’s inequality that
\begin{equation}
\label{Mtconv}
\Pp\big(|b_n^{-1}Q_n(b_nt)-\alpha t|>\veps) \le \frac{b_n^{-2}\mathrm{Var}(Q_n(b_nt))}{\veps^2}\xrightarrow{n\to\infty}0.
\end{equation}
The uniform convergence now follows from the monotonicity of $Q_n$ and can be deduced in a similar way as in the previous proof. 
\end{proof}

\begin{proof}[Proof of \cref{LLN lemma}]
We apply Lemma~\ref{lem: LLN-gen} with $x_n(i)\equiv 1$ for all $i\in [n]$, $y_n(i)=(1+\lambda n^{-1/3})\wt{i}/|\wtv|_1$ and $b_n=n^{2/3}$. Then both conditions in~\eqref{cond: LLN-gen} are fulfilled with $\alpha=1$, due to the assumption~\eqref{maximal weight};~\eqref{cond: bn} also holds true as $\mathbf x_n$ is constant in this case.  
\end{proof}

\begin{proof}[Proof of Lemma~\ref{lem:lln_roots}]
The first two convergences can be deduced by applying Lemma~\ref{lem: LLN-gen}. Indeed, for~\eqref{Vn conv}, we take $x_n(i)=\wt{i}$, $y_n(i)=\wt{i}/|\wtv|_1$ and $b_n=n^{1/3}$. Then both limits in~\eqref{cond: LLN-gen} hold true with $\alpha=1$, thanks to~\eqref{first moment},~\eqref{second moment} and~\eqref{maximal weight}. Moreover, in this case,~\eqref{cond: bn} is simply~\eqref{maximal weight}. An application of Lemma~\ref{lem: LLN-gen} yields~\eqref{Vn conv}. Similarly for~\eqref{Vn conv 2}, we take $x_n(i)=\wt{i}$, $b_n=n^{1/3}$ as previously; let also $y_n(i)=1/n$. Now the first convergence in~\eqref{cond: LLN-gen} with a choice of $\alpha=\mu$ follows from~\eqref{first moment} and the second one follows from the definition of $\mathbf y_n$. Moreover, ~\eqref{cond: bn} holds for precisely the same reason as before. We obtain~\eqref{Vn conv 2} from the lemma. 

For the two remaining convergences, we introduce
\[
A(t)=\sum_{i\in [n]}\I{E_i\le t} \quad \text{and}\quad \hat A(t)=\sum_{i\in [n]}\I{\hat E_i\le t}, \quad t\ge 0. 
\]
We note that $T_k=\inf\{t\ge 0: A(t)\ge k\}$ and $\hat T_k=\inf\{t\ge 0: \hat A(t)\ge k\}$, for each $k\ge 1$. 
Another application of Lemma~\ref{lem: LLN-gen} yields that $(n^{-1/3}A(n^{1/3}t), t\ge 0)$ converges in probability to $(t, t\ge 0)$ uniformly on compacts. 
Since the limit process is continuous, strictly increasing,~\eqref{Tn conv} holds true. The convergence in~\eqref{Tn conv2} can be argued similarly. 
\end{proof}

\begin{claim}\label{claim:inf_W}
    For any $\veps>0$ and any $C>0$, 
    \[ \limsup_{T\to\infty}\Pp\left(\inf_{0\le t\le T/2} \bm(t)-\inf_{0\le t\le 2T} \bm(t)\le C  \right)<\veps.\]
\end{claim}
\begin{proof}
In this proof, we write $\alpha=\sqrt{\tfrac{\mu'}{\mu}}$ and $f(t)=\tfrac{1}{2}\tfrac{\mu'}{\mu^2}t^2 -\lambda t$ so that $\bm(t)=\alpha B(t)-f(t)$ for $B(t)$ a standard Brownian motion. 

As $f$ is strictly convex, it follows that for all $T$ large enough,  $\sup_{0\le t \le T/2}f(t)=f(T/2)$, and note that for any $x>0$, 
\begin{align*}
&\Pp\left(\inf_{0\le t\le T/2} \bm(t)-\inf_{0\le t\le 2T} \bm(t)\le C  \right)\le \Pp\left(\inf_{0\le t\le T/2} \bm(t) \le -x\right) + \Pp(\bm(2T) \ge -x-C)\\
&\quad\le \Pp\left(\inf_{0\le t\le T/2} \alpha B(t) \le f(T/2)-x\right)+\Pp\left(\alpha B(2T) \ge f(2T)-x-C\right)\\
&\quad=2\Pp(\alpha B(T/2) \le f(T/2)-x)+\Pp(\alpha B(2T) \ge f(2T)-x-C),
\end{align*}
where the last line follows from the reflection principle for Brownian motion. Setting $x=2f(T/2)$, we see that the last line equals 
\[2\Pp\left(N \le \frac{-f(T/2)}{\alpha \sqrt{\tfrac{T}{2}}}\right)+\Pp\left(N \ge \frac{f(2T)-2f(T/2)-C}{\alpha \sqrt{2T}}\right), \]
where $N$ is a standard normal random variable. Finally, it is easy to see from the definition of $f(t)$ that $f(T/2)/(\alpha\sqrt{T/2})$ and $(f(2T)-2f(T/2)-C)/(\alpha \sqrt{2T})$ both tend to $\infty$ as $T\to \infty$, which together with the preceding bound implies the claim. 
\end{proof}

\bibliographystyle{plainnat}
\bibliography{template}

\end{document}